\documentclass[11pt,reqno]{amsart}

\usepackage{amsmath, amsthm, amssymb}
\usepackage{amsfonts}

\usepackage{hyperref}

\usepackage[ansinew]{inputenc}
\usepackage[dvips]{epsfig}
\usepackage{graphicx}
\usepackage[english]{babel}

\usepackage{thmtools}
\theoremstyle{plain}
\declaretheorem[title=Theorem, parent=section]{theorem}
\declaretheorem[title=Lemma,sibling=theorem]{lemma}
\declaretheorem[title=Proposition,sibling=theorem]{proposition}

\theoremstyle{definition}
\declaretheorem[title=Definition,sibling=theorem]{definition}
\declaretheorem[title=Remark,sibling=theorem]{remark}
\declaretheorem[title=Remark, numbered=no]{remark*}

\declaretheorem[title=Assumption, numbered=no]{assumption*}

\numberwithin{equation}{section}

\usepackage[backgroundcolor=white, bordercolor=blue,
linecolor=blue]{todonotes}

\usepackage{dsfont}
\usepackage{bbm}

\newcommand{\N}{\mathbb{N}}
\newcommand{\R}{\mathbb{R}}

\newcommand{\cE}{\mathcal{E}}

\newcommand{\cM}{\mathcal{M}}

\newcommand{\al}{\alpha}

\newcommand{\dl}{\delta}

\newcommand{\lm}{\lambda}
\newcommand{\Lm}{\Lambda}

\newcommand{\varep}{\varepsilon}
\newcommand{\vp}{\varphi}
\newcommand{\sig}{\sigma}

\newcommand{\om}{\omega}
\newcommand{\Om}{\Omega}
\newcommand{\z}{\zeta}

\newcommand{\eps}{\varepsilon}

\newcommand{\loc}{\mathrm{loc}}

\DeclareMathOperator{\supp}{supp}

\DeclareMathOperator{\dvg}{div}

\DeclareMathOperator*{\osc}{osc}

\DeclareMathOperator{\tail}{Tail}

\renewcommand{\d}{\textnormal{\,d}}

\newcommand{\average}{{\mathchoice {\kern1ex\vcenter{\hrule height.4pt
width 6pt depth0pt} \kern-9.7pt} {\kern1ex\vcenter{\hrule
height.4pt width 4.3pt depth0pt} \kern-7pt} {} {} }}
\newcommand{\dashint}{\average\int}

\def\Xint#1{\mathchoice
    {\XXint\displaystyle\textstyle{#1}}%
    {\XXint\textstyle\scriptstyle{#1}}%
    {\XXint\scriptstyle\scriptscriptstyle{#1}}%
    {\XXint\scriptscriptstyle\scriptscriptstyle{#1}}%
    \!\int}
\def\XXint#1#2#3{\setbox0=\hbox{$#1{#2#3}{\int}$}
    \vcenter{\hbox{$#2#3$}}\kern-0.5\wd0}
\def\bint{\Xint-}
\def\dashint{\Xint{\raise4pt\hbox to7pt{\hrulefill}}}
\def\tmint{\Xint{\raise0pt\hbox to6pt{\hrulefill}}}

\def\XXiint#1#2#3{\setbox0=\hbox{$#1{#2#3}{\iint}$}
    \vcenter{\hbox{$#2#3$}}\kern-0.5\wd0}

\begin{document}
\allowdisplaybreaks
\title[Nonnegative solutions to nonlocal parabolic equations]{Nonnegative solutions to nonlocal parabolic equations}
 
\author{Naian Liao}
 
\author{Marvin Weidner}

\address{Fachbereich Mathematik, Paris-Lodron-Universit\"at Salzburg, Hellbrunner Str, 24, 5020 Salzburg, Austria}
\email{naian.liao@plus.ac.at}

\address{Departament de Matem\`atiques i Inform\`atica, Universitat de Barcelona, Gran Via de les Corts Catalanes 585, 08007 Barcelona, Spain}
\email{mweidner@ub.edu}
\urladdr{https://sites.google.com/view/marvinweidner/}



\allowdisplaybreaks

\begin{abstract}
We aim to study nonnegative, global solutions to a general class of nonlocal parabolic equations with bounded measurable coefficients. First, we prove a Widder-type theorem. Such a result has previously been studied only for certain translation invariant operators, and new ideas are needed in our general setting. Second, we establish sharp two-sided bounds for the fundamental solution via purely variational techniques, entirely bypassing tools from semigroup theory, Dirichlet forms, and stochastic analysis. Third, we derive sharp Harnack-type estimates that are novel even for the fractional heat equation. 
\end{abstract}

\allowdisplaybreaks

\maketitle


\section{Introduction}
Let us consider the following global parabolic equation
\begin{align}
\label{eq:local-PDE}
    \partial_t u - \dvg(\mathbf{A}(t,x)\nabla u) = 0 ~~ \text{ in } (0,T) \times \R^d,
\end{align}
where $\mathbf{A} : (0,T) \times \R^d \to \R^{d \times d}$ is a symmetric matrix which is bounded, measurable, and uniformly elliptic. A celebrated result by D.G. Aronson \cite{Aro68,Aro71} establishes the existence of the \emph{fundamental solution} $\Gamma_t(x,y) : (0,T) \times \R^d \times \R^d \to [0,\infty]$ to \eqref{eq:local-PDE} with the following properties:
\begin{itemize}
    \item[(i)] For any nonnegative weak solution $u$ to \eqref{eq:local-PDE}, there is a unique Radon measure $\mu$ such that $u$ satisfies $u(0) = \mu$ and enjoys the following representation:
    \begin{align}
    \label{eq:local-rep}
        u(t,x) = \int_{\R^d} \Gamma_t(x,y) \d\mu(y).
    \end{align}
    \item[(ii)] For any Radon measure $\mu$ with proper growth at infinity, the function $u$ defined by \eqref{eq:local-rep} satisfies $u(0) = \mu$ and is a weak solution to \eqref{eq:local-PDE}.
    \item[(iii)] The fundamental solution $\Gamma_t(x,y)$ satisfies the following two-sided Gaussian bounds:
    \begin{align}
        \label{eq:Gaussian}
        C_1 t^{-\frac{d}{2}} e^{-c_1\frac{|x-y|^2}{t} } \le \Gamma_t(x,y) \le C_2 t^{-\frac{d}{2}} e^{-c_2\frac{|x-y|^2}{t} } \qquad \forall\, t \in (0,T), ~~ x,y \in \R^d,
    \end{align}
    where $C_1,C_2,c_1,c_2$ depend only on $d$ and the ellipticity constants of $\mathbf{A}$.
\end{itemize}

This result is of central importance in the theory of linear, second order parabolic equations as it establishes well-posedness of the Cauchy problem and shows that all nonnegative weak solutions to \eqref{eq:local-PDE} are given by the ``convolution" of its initial data $\mu$ with a kernel $\Gamma_t(x,y)$. Prior to \cite{Aro68,Aro71}, (i) and (ii) had been studied for the classical heat equation in 1D by Widder \cite{Wid44}, and for equations in non-divergence form by Krzy\.za\'nski \cite{Krz64}, excluding however solutions with measure initial data. The representation formula \eqref{eq:local-rep} is a deep structural result for \eqref{eq:local-PDE}, effectively reducing its analysis to the study of the fundamental solution. The two-sided estimates \eqref{eq:Gaussian} establish that $\Gamma_t(x,y)$ is comparable to the fundamental solution of the classical heat equation in a pointwise way, up to multiplicative constants. In this sense, the bounds \eqref{eq:Gaussian} are robust within the class of second order divergence form operators that have bounded, measurable, and uniformly elliptic coefficients. 
Needless to say, the results (i), (ii), and (iii) are also vital in other areas of mathematics, such as geometric analysis and probability theory, where $\Gamma_t(x,y)$ arises, for instance, as the heat kernel on a manifold or as the transition density of a stochastic diffusion process, respectively.

\subsection{Nonlocal equations with bounded measurable coefficients}

The goal of this article is to establish an analogous theory -- including (i), (ii), and (iii) -- for \emph{nonlocal} parabolic equations with bounded, measurable, and uniformly elliptic coefficients.

A nonlocal counterpart of the global parabolic equation \eqref{eq:local-PDE} is given by  
\begin{equation}\label{Eq:1:1}
    \partial_t u(t,x)-\mathcal{L}_t u(t,x)=0\quad\text{in}\quad (0,T) \times \R^d,
\end{equation}
where $\mathcal{L}_t$ is an integro-differential operator of the form
\begin{align*}
    -\mathcal{L}_t u(t, x) = \text{p.v.}~\int_{\R^d} (u(t, x) - u(t, y)) K(t;x,y) \d y
\end{align*}
for a measurable kernel $K: (0,T) \times \R^d \times \R^d \to [0,\infty]$ that is symmetric in $(x,y)$. We are interested in operators $\mathcal{L}_t$ that are modeled upon the fractional Laplacian $-(-\Delta)^s$, $s \in (0,1)$. In this context, natural uniform ellipticity conditions on $K$, given $0 < \lambda \le \Lambda < \infty$, are the following upper bound
\begin{equation}\label{eq:up-bd}
    K(t;x,y)  \le \frac{\Lambda}{|x-y|^{d+2s}}
\end{equation}
for any $t\in (0,T)$, and $x,y \in \R^d$, and the following coercivity condition
\begin{equation}\label{eq:coercive}
    \iint_{B_r(x_o)\times B_r(x_o)} |v(x) - v(y)|^2K(t;x,y) \d y \d x \ge \lm\, [v]^2_{H^s(B_r(x_o))}
\end{equation}
for any $t\in (0,T)$, $r > 0$, $x_o \in \R^d$ and $v\in H^s(B_r(x_o))$.  

The coercivity condition \eqref{eq:coercive} is a substantially weaker assumption than a pointwise lower bound of the kernel, namely, that for any $t\in (0,T)$, $x,\,y\in\R^d$ it holds
\begin{equation}\label{eq:low-bd}
    K(t;x,y) \ge \frac{\lambda}{|x-y|^{d+2s}}.
\end{equation}

In the last twenty years there has been a huge interest in parabolic equations governed by nonlocal diffusion operators with measurable coefficients, such as \eqref{Eq:1:1}. A cornerstone in this area is the celebrated De Giorgi-Nash-Moser theory, whose nonlocal variant has been successfully developed in \cite{Kas09,DKP14,DKP16,Coz17,CKW22b,CKW23,BKO23} for elliptic problems and in \cite{CCV11,FeKa13,KaWe23b,APT22,Lia24,Lia24b,Lia24c} for parabolic problems.

Despite recent progress, a Widder-type theory as in (i) and (ii), as well as two-sided estimates for the fundamental solution in analogy to (iii), remain open for nonlocal equations \eqref{Eq:1:1} in the general setting presented above. So far, such results have only been obtained in certain special cases, as we will discuss in more detail below. In this article, we develop a unified theory for nonlocal equations \eqref{Eq:1:1} under the general ellipticity assumptions \eqref{eq:up-bd}, \eqref{eq:coercive}, and establish complete nonlocal analogs of (i), (ii), and (iii). A key strength of our approach lies in its foundation on purely variational tools.

\subsection{A nonlocal Widder-type theory}

Our first main result proves a Widder-type theorem for nonlocal equations \eqref{Eq:1:1} with bounded measurable coefficients satisfying \eqref{eq:up-bd}, \eqref{eq:coercive}. Under these general assumptions on $K$ we show the existence of the fundamental solution $(t,x) \mapsto p_t(x,y)$ (see Subsection \ref{subsec:ex-fs}), which solves for any $y \in \R^d$
\begin{equation*}
\left\{\begin{array}{rcl}
\partial_t p_t(\cdot,y) -\mathcal{L}_t p_t(\cdot,y) &=& 0 ~~~~~~ \text{ in }  (0,T) \times \R^d, \\
p_t(\cdot,y) &=& \delta_{\{y\}} ~~ \text{ on } \{ 0 \} \times \R^d.
\end{array}\right.
\end{equation*}

With the fundamental solution at hand, we are able to prove the following Widder-type theorem in analogy to (i) and (ii).

\begin{theorem}[Widder-type theorem]
    \label{thm:widder}
    Assume that the kernel $K$ satisfies the upper bound \eqref{eq:up-bd} and the coercivity condition \eqref{eq:coercive}. Then, the following hold true:
    \begin{itemize}
        \item[(i)] Let $u$ be a nonnegative, global weak solution to 
        \begin{align}
        \label{eq:PDE-Widder}
            \partial_t u(t,x)-\mathcal{L}_t u(t,x)=0\quad\text{in}\quad (\eps,T) \times \R^d
        \end{align}
        for any $\varep\in(0,T)$ in the sense of \autoref{def:global-sol}. Then, there exists a unique, nonnegative Radon measure $\mu$ with 
    \begin{equation}\label{eq:mu-growth}
            \int_{\R^d} \frac{ \d |\mu|(x)}{(1 +|x|)^{d+2s}} < \infty
    \end{equation}
     such that $u(0) = \mu$ in the sense of vague convergence of measures, i.e.
     \begin{align}
     \label{eq:init-intro}
         \int_{\R^d} \psi(x) \d \mu(x) = \lim_{t \searrow 0} \int_{\R^d} \psi(x) u(t,x) \d x, ~~ \forall\, \psi \in C_c(\R^d),
     \end{align}
     and it holds
\begin{align}\label{eq:u-rep}
    u(t,x) = \int_{\R^d} p_t(x,y) \d \mu(y).
\end{align}
\item[(ii)] The representation \eqref{eq:u-rep} with any -- possibly signed -- Radon measure $\mu$ satisfying \eqref{eq:mu-growth} defines a global weak solution to \eqref{eq:PDE-Widder} for any $\eps \in (0,T)$ in the sense of \autoref{def:global-sol} and it holds $u(0) = \mu$ in the sense of \eqref{eq:init-intro}.
    \end{itemize}
\end{theorem}

A key feature of (i) is the identification of a uniquely defined initial trace $\mu$ for nonnegative weak solutions. Let us highlight that we can recover such an initial condition even for solutions that are \emph{local in time}, meaning that no regularity is assumed up to the initial time $t=0$. In fact, we only require solutions to be in the following natural function space:
\begin{align}
\label{eq:local-in-time}
    u \in L^{\infty}_{\loc}((0,T);L^2_{\loc}(\R^d)) \cap L^2_{\loc}((0,T);H^s_{\loc}(\R^d)) \cap L^1_{\loc}((0,T);L^1(\R^d;w)).
\end{align}
Here, $L^1(\R^d;w)$ denotes the natural weighted $L^1$ space (see \eqref{eq:L1w}) where $w(x) := (1 + |x|)^{-d-2s}$ matches the decay of $K$ at infinity. We refer to \autoref{def:global-sol} for a precise definition of the solution concept which is present in \autoref{thm:widder}.

\begin{remark}
    Note that we are not restricted to nonnegative solutions in \autoref{thm:widder}(ii). It is natural to ask whether also \autoref{thm:widder}(i) holds true for \emph{signed weak solutions}. In \autoref{lemma:uniqueness-signed} we prove that signed weak solutions to \eqref{Eq:1:1} with a given initial trace are unique in the function space
    \begin{align*}
            L^{\infty}((0,T);L^2_{\loc}(\R^d)) \cap L^2((0,T);H^s_{\loc}(\R^d)) \cap L^1((0,T);L^1(\R^d;w)).
    \end{align*}
    This space, as opposed to \eqref{eq:local-in-time}, prescribes some regularity up to time $t=0$. Solutions can be constructed in such space, provided that sufficiently regular initial data are prescribed (see \cite[Appendix~A]{LiWe24}). It is an interesting open problem whether signed solutions with a given measure initial trace are unique in the larger space \eqref{eq:local-in-time}, even in the simplest case $\mathcal{L}_t = -(-\Delta)^s$. 
\end{remark}

Let us now discuss to what extend \autoref{thm:widder} is new.
\begin{itemize}
    \item A Widder-type theory has first been established in the special case $\mathcal{L}_t = -(-\Delta)^s$ in \cite{BPSV14} and has been extended to measure initial data in \cite{BSV17}. In these articles, the authors analyze well-posedness of \eqref{Eq:1:1} for various notions of weak, distributional, and strong solutions. All of their results are stated under stronger regularity assumptions at the initial time, namely, $u \in L^1((0,T);L^1(\R^d;w))$, even for nonnegative solutions.
    
    \item In the recent preprint \cite{GQSV25}, the authors study a Widder-type theory for a class of translation invariant kernels $K(t;x,y) = K(x-y)$ satisfying mixed polynomial growth conditions from above and below. These bounds include operators of variable order which are not treated in our article. At the same time, their setting does not allow for integrated lower bounds as in \eqref{eq:coercive}, and therefore the kernels discussed in \autoref{rem:gen-lb} fall outside their scope. They do not assume any regularity at the initial time in case solutions are nonnegative, similar to \autoref{thm:widder}(i).
\end{itemize}

Our result \autoref{thm:widder} is the first to establish a Widder-type theory for nonlocal equations with bounded measurable coefficients, as all previous results are restricted to translation invariant nonlocal operators.

Since we do not impose any regularity assumption on the kernel $K$, the equation \eqref{eq:PDE-Widder} has to be understood in the weak sense. Hence, we are restricted to the use of variational techniques, which is in stark contrast to the recent works \cite{BPSV14,BSV17,GQSV25}. In fact, their approaches essentially use that the operators $\mathcal{L}$ can be evaluated in a classical way for smooth functions (see \cite[Lemma 4.2, Lemma 11.1]{BSV17} and \cite[Lemma 3.7, Lemma 3.9]{GQSV25}), which allows them to use barrier arguments.

In our purely variational framework, barrier arguments are no longer available. Therefore, several new ideas are needed. Most importantly, we employ fine local estimates, such as the improved weak Harnack inequality (see \autoref{thm:improved-weak-Harnack}), which we proved in a previous work \cite{LiWe24}. This allows us to propagate information on the $L^1$-norms of solutions forward time. A delicate analysis is required to maintain control over the global behavior of the solution.

In \cite[Proposition A.1]{LiWe24} we have already shown that \eqref{eq:u-rep} is the unique solution to \eqref{Eq:1:1} with a given initial datum in $L^2(\R^d)$. This setting is in alignment with the theory of $L^2$-semigroups where solutions are more regular up to the initial time and attain their initial datum in $L^2(\R^d)$. A key contribution of the current paper is to lift the theory from $L^2(\R^d)$ to measure initial data merely satisfying \eqref{eq:mu-growth}, a setting well beyond the reach of classical semigroup methods.

Let us close this part of the introduction by mentioning several related works. Widder-type existence and uniqueness results have also been studied for nonlinear

\subsection{Two-sided estimates for the fundamental solution}

Our second main result establishes that the fundamental solution $p_t(x,y)$ to \eqref{Eq:1:1} enjoys two-sided polynomial bounds, whenever $K$ satisfies the pointwise upper and lower bounds \eqref{eq:up-bd} and \eqref{eq:low-bd}. This result is a complete analog of (iii) in the local case.

\begin{theorem}
    \label{thm:hkb}
    Assume that the kernel $K$ satisfies the upper bound \eqref{eq:up-bd} and the lower bound \eqref{eq:low-bd}. Then, the fundamental solution $p_t(x,y)$ satisfies the following two-sided bounds:
    There exist $c_1,c_2 > 0$, depending only on the data $\{d,s,\lambda,\Lambda\}$, such that 
\begin{align}
\label{eq:hkb}
c_1 \left( t^{-\frac{d}{2s}} \wedge \frac{t}{|x-y|^{d+2s}}  \right) \le p_t(x,y) \le  c_2 \left( t^{-\frac{d}{2s}} \wedge \frac{t}{|x-y|^{d+2s}}  \right)
\end{align}
holds for any $x,\,y\in\R^d$ and $t \in (0,T)$. Moreover, the upper bound in \eqref{eq:hkb} remains true if $K$ merely satisfies \eqref{eq:coercive} instead of \eqref{eq:low-bd}.
\end{theorem}

Estimate \eqref{eq:hkb} states that the fundamental solution to any nonlocal operator with bounded measurable coefficients is comparable to the fundamental solution to the fractional heat equation $\partial_t u + (-\Delta)^s u = 0$, see \cite{BlGe60}. It is noteworthy that the upper bound in \eqref{eq:hkb} continues to hold under the weaker coercivity condition \eqref{eq:coercive}, whereas the lower bound might break if \eqref{eq:low-bd} is violated; see Subsections~\ref{subsec-upper} \& \ref{subsec-lower} for a more detailed discussion.

There exists a wide literature on estimates for the fundamental solution to nonlocal parabolic problems already in very general settings that include operators with bounded measurable coefficients, as we discuss below. The main new contributions of our article are twofold:
\begin{itemize}
    \item In contrast to existing results, we establish \eqref{eq:hkb} for operators whose kernels are allowed to be time-dependent.
    
    \item While previous approaches are based on powerful tools from semigroup theory, Dirichlet forms, or stochastic analysis, our technique is entirely variational. Our proof of the lower bound in \eqref{eq:hkb} is completely new and considerably shorter than all previous ones.
\end{itemize}

 We briefly review the development of fundamental solution estimates for nonlocal operators. Two-sided bounds of the type \eqref{eq:hkb} were obtained in \cite{BaLe02} and \cite{ChKu03} within a probabilistic framework, under the assumption that $K$ is independent of time and satisfies \eqref{eq:up-bd} and \eqref{eq:low-bd}. Their derivation of the upper estimate relies crucially on the extension of Davies' method to jump processes in \cite{CKS87}, which is in turn based on results from the theory of Markov semigroups. The lower estimate heavily uses the L\'evy system formula, a powerful tool from stochastic analysis. Building on these ideas, \cite{ChKu08}, \cite{BGK09}, and \cite{CKW21} established two-sided bounds as in \eqref{eq:hkb} under more general assumptions on $K$, allowing for kernels with mixed polynomial growth, and replacing $\R^d$ by doubling metric measure spaces.

Further generalizations were developed in a series of papers including \cite{GHL14}, \cite{GHH17}, and \cite{GHH18}, where the theory was extended to metric measure spaces with walk dimension greater than two. These works provide characterizations of upper and two-sided heat kernel bounds via geometric and analytic conditions on the kernel and the underlying space. Notably, their method avoids probabilistic tools and instead uses comparison principles for heat semigroups and Dirichlet form theory. 

We also refer to \cite{CKKW21,CKW22c,GHH23,GHH24} for more recent results on heat kernel estimates in even more general settings where the upper and lower bounds in \eqref{eq:up-bd} and \eqref{eq:low-bd} contain scaling functions that can vary highly in space. Their proofs use variations of the aforementioned techniques. 

Note that most of these approaches do not seem to generalize to time-dependent kernels in a straightforward way. In fact, they rely partly on regularity results of De Giorgi-Nash-Moser type for \emph{elliptic nonlocal equations} which are employed in a parabolic context by using the relation between the heat semigroup and the resolvent operator. These connections are only available if kernels are independent of time (see for instance \cite[Lemma 5.6]{GHH18} and \cite[Lemma 4.14, 4.15]{CKW21}), which we are able to drop completely in this article. Alternatively, \cite{MaMi13a,MaMi13b} provide the upper heat kernel estimate in \eqref{eq:hkb} for time-dependent kernels, assuming the more restrictive pointwise lower bound \eqref{eq:low-bd} instead of \eqref{eq:coercive}. These works focus on the treatment of an additional first-order drift term that is divergence-free.

As was already mentioned, a key strength of our approach is that it entirely bypasses tools from semigroup theory, Dirichlet forms, and stochastic analysis. This is in stark contrast to all articles that were mentioned above. We highlight that our proof of the lower bound in \eqref{eq:hkb} seems to be completely new and is considerably shorter than all previous ones. It merely relies on the nonlocal parabolic weak Harnack inequality with a tail term that was established in \cite{KaWe23b} (see also \cite{Lia24}). We believe this technique to be applicable in more general contexts and, as an instance, we demonstrate in \autoref{rem:gen-lb} that it carries over to kernels $K$ violating the pointwise lower bound \eqref{eq:low-bd}.

To prove the upper estimate in \eqref{eq:hkb}, we apply the results in \cite{KaWe23}, where Aronson's technique \cite{Aro68} was extended to nonlocal operators with kernels satisfying \eqref{eq:up-bd} and \eqref{eq:coercive}. At the time \cite{KaWe23} was written, the fundamental solution $p_t(x,y)$ had not been  constructed for time-dependent problems, and therefore, some of its properties were posed as assumptions in the main theorem of \cite{KaWe23}. A new contribution of the current article is to construct the fundamental solution and to rigorously verify all its required properties (see Subsection \ref{subsec:ex-fs}).

\subsection{Harnack-type estimates for global solutions}

Time-insensitive Harnack estimates have been established in \cite{LiWe24} within the same framework as in this article. It is interesting to observe that they also follow from the Widder-type theorem and the two-sided bounds of the fundamental solution we prove in the current work. Next, we go beyond and show sharp estimates on the growth of nonnegative global solutions to \eqref{Eq:1:1} as $|x| \to \infty$. In particular, we have the following elliptic-type Harnack inequality \eqref{eq:Harnack-type-estimate} that explicitly estimates the quotient of two values of solutions, generalizing the result in \cite{LiWe24}.

\begin{theorem}[Harnack-type estimate]
\label{thm:est-global}
        Assume that the kernel $K$ satisfies the upper bound \eqref{eq:up-bd} and the lower bound \eqref{eq:low-bd}. Let $u$ be a nonnegative global solution to \eqref{eq:PDE-Widder} for any $\varep\in(0,T)$ in the sense of \autoref{def:global-sol}. Then, there exists $c>0$ depending only on the data $\{d,s,\lambda,\Lambda\}$, such that the estimate
        \begin{align}
        \label{eq:Harnack-type-estimate}
            \frac{u(t,x)}{u(\tau,y)}\le c\,\left(\frac{t}{\tau}\right)^{-d/2s}\bigg(1+\frac{|x-y|}{\tau^{1/2s}}\bigg)^{d+2s}
        \end{align}
        holds true for any $0<t \le \tau < T$ and $x,\,y\in \R^d$.
\end{theorem}

\begin{remark}\label{rmk:sub-potential}
The Harnack estimate \eqref{eq:Harnack-type-estimate} automatically yields a sub-potential lower bound for any nonnegative global solution $u$, namely,
\[
\frac{u(t,x)}{u(t,0)}\ge c \, \displaystyle\bigg(1+\frac{|x|}{t^{1/2s}}\bigg)^{-(d+2s)}
\]
for any $x\in\R^d$ and $t>0$.
\end{remark}

\begin{remark}The Harnack-type estimate \eqref{eq:Harnack-type-estimate} is sharp as testified by the fundamental solution. Indeed,
consider the function
\[
P(t,x):=\frac{t^{-d/2s}}{\displaystyle\bigg(1+\frac{|x|}{t^{1/2s}}\bigg)^{d+2s}}
\]
defined for $x\in\R^d$ and $t>0$. It is straightforward to check that $P(t,x)$ is comparable to $p_t(x,0)$ thanks to the two-sided bounds in \eqref{eq:hkb}.
On the other hand, it is obvious that
\[
\frac{P(t,0)}{P(t,x)}=\displaystyle\bigg(1+\frac{|x|}{t^{1/2s}}\bigg)^{d+2s}.
\]
Hence, applying the two-sided bounds in \eqref{eq:hkb}, the sharpness of \eqref{eq:Harnack-type-estimate} is testified by the fundamental solution $p_t(x,0)$.
\end{remark}

Previously, a sharp Harnack-type estimate was only studied for the fractional heat equation, see \cite{BSV17, DeSi23,WeZa23}. However, existing results rely on Widder-type theorems that exclude the fundamental solution. Strikingly, our study demonstrates the following: 
\begin{itemize}
    \item A sharp Harnack-type estimate holds true for any nonnegative global weak solution in the general class \eqref{eq:local-in-time}, regardless of its initial datum. In particular, our estimates are new even for the fractional heat equation, and they can be proved without any Widder-type theorem.

    \item Such an estimate is a structural property of nonnegative global solutions, paralleling a well-known result of Moser (\cite[(1.7)]{Mos64}). 
    \item The time-insensitive Harnack estimates established in \cite{LiWe24} are sufficient to describe the local and global behavior of nonnegative solutions within the framework of De Giorgi-Nash-Moser theory.
\end{itemize}

\begin{remark}
    It would be interesting to know if one can replace \eqref{eq:low-bd} by the weaker condition \eqref{eq:coercive} and what the sharp Harnack estimate will be in that case. One possible strategy to achieve such generalization could be to extend the techniques in \cite{LiWe24}. Alternatively, since our Widder-type theory already includes the general class of kernels satisfying merely the coercivity assumption \eqref{eq:coercive}, an extension of \autoref{thm:est-global} could also be achieved by combining \autoref{thm:widder} with two-sided heat kernel estimates for more general nonlocal operators. See \autoref{rem:gen-lb} for further discussion.
\end{remark}

As an application of our Widder-type theorem (see \autoref{thm:widder}), we can also estimate the growth of solutions in terms of their initial datum.

\begin{theorem}[Growth estimate via initial data]
\label{thm:est-global-}
        Assume that the kernel $K$ satisfies the upper bound \eqref{eq:up-bd} and the coercivity condition \eqref{eq:coercive}. Let $u$ be a nonnegative global solution to \eqref{eq:PDE-Widder} for any $\varep\in(0,T)$ in the sense of \autoref{def:global-sol} and let $\mu$ be the initial measure of $u$ identified in \autoref{thm:widder}. Then, there exists $c>0$ depending only on the data $\{d,s,\lambda,\Lambda\}$, such that
        \begin{align}\label{eq:Harnack-initial}
            u(t,x) \le c  t^{-d/2s} \int_{\R^d} \frac{(t^{1/2s} +|x|)^{d+2s}}{(t^{1/2s} + |y|)^{d+2s}}\d \mu(y).
        \end{align}
        holds true for any $t\in(0,T]$ and $x \in \R^d$.
\end{theorem}

\subsection{Outline}
This article is structured as follows. In Section \ref{sec:prelim} we introduce some notation and define the weak solution concepts present in this article. In Section \ref{sec:wHI}, we recall several versions of weak Harnack inequalities for nonnegative supersolutions to nonlocal equations. Section \ref{sec:fund-sol} is dedicated to the construction of the fundamental solution and the proof of its two-sided bounds (see \autoref{thm:hkb}). In Section \ref{sec:Widder} we establish the Widder-type theorem for nonlocal equations with bounded measurable coefficients (see \autoref{thm:widder}). In Section \ref{sec:est-global}, we prove \autoref{thm:est-global} and \autoref{thm:est-global-}. Finally, this paper contains two appendices. In Appendix \ref{sec:signed-sol} we prove a uniqueness result for signed solutions in a special class (see \autoref{lemma:uniqueness-signed}). Finally, Appendix \ref{sec:cont-initial} contains a proof of continuity of solutions up to the initial level given some continuous initial data (see \autoref{prop:initial}).

\subsection*{Acknowledgments}

Naian Liao was supported by the FWF-project P36272-N ``On the Stefan type problems".
Marvin Weidner has received funding from the European Research Council (ERC) under the Grant Agreement No 801867 and the Grant Agreement No 101123223 (SSNSD), and by the AEI project PID2021-125021NA-I00 (Spain).

\section{Preliminaries}
\label{sec:prelim}

In this article, $B_R(x_o)$ denotes the ball of radius $R$ and center $x_o$ in $\R^d$, whereas $K_R(x_o)$ denotes the cube of side length $2R$, center $x_o$ and faces parallel with the coordinate plane. We refer to the weight $w(x):=(1+|x|)^{-d-2s}$ and define the weighted Lebesgue space
\begin{align}
\label{eq:L1w}
    L^1(\R^d;w) = \bigg\{ u \in L^1_{\loc}(\R^d) : \Vert u \Vert_{L^1(\R^d;w)} = \int_{\R^d} \frac{|u(x)|}{(1+|x|)^{d+2s}}\d x < \infty \bigg\}.
\end{align}
Moreover, we define a bilinear form
\begin{align*}
\cE^{(t)}(u,v) := \int_{\R^d} \int_{\R^d} (u(x) - u(y))(v(x)-v(y))K(t;x,y) \d y \d x
\end{align*}
which will be finite if $u$ and $v$ possess certain regularity.

Let us introduce the solution concepts. The first one is a local notion.

 \begin{definition}[local solution]\label{def:local-sol}
 Let $I \subset \R$ be an interval and $\Omega \subset \R^d$ be a bounded domain.
We say that 
$$
u \in L^{\infty}(I;L^2_{\loc}(\Omega)) \cap L^2(I;H^s_{\loc}(\Omega)) \cap L^1(I;L^1(\R^d;w))
$$ is a weak supersolution to
\begin{align}
\label{eq:PDE-domains}
    \partial_t u - \mathcal{L}_t u = 0 ~~ \text{ in } I \times \Omega,
\end{align}
if for any $\phi \in H^1(I;L^2(\Om))\cap  L^2(I;H^s(\R^d))$ with $\phi \ge 0$ and $\supp(\phi) \subset I \times \Omega$ it holds
\begin{align*}
        -\int_I \int_{\R^d} u(t,x) \partial_t \phi(t,x) \d x \d t + \int_I \cE^{(t)}(u(t),\phi(t)) \d t \ge 0.
    \end{align*}
    We say that $u$ is a weak subsolution if $-u$ is a weak supersolution, whereas $u$ is a weak solution if $u$ is a weak supersolution and a weak subsolution.
\end{definition}

The second notion is a global one. Yet, we do not require solutions to be in  $H^s(\R^d)$ nor in $L^2(\R^d)$; otherwise the notion would exclude solutions that grow at infinity.
\begin{definition}[global solution]\label{def:global-sol}
    Let $I$ be an interval. We say that
    $$
    u \in L^{\infty}(I;L^2_{\loc}(\R^d)) \cap L^2(I;H^s_{\loc}(\R^d)) \cap L^1(I;L^1(\R^d;w))
    $$ is a global weak solution to 
    \begin{align*}
    \partial_t u - \mathcal{L}_t u = 0 ~~ \text{ in } I \times \R^d,
\end{align*}
if it is a local weak solution in $I\times B_R$ for any $R>0$.
\end{definition}

By $\cM$ we denote the set of signed Radon measures $\mu$ on $\R^d$. Meanwhile, we say that $\mu \in \cM_+$ if $\mu$ is nonnegative.

Next, we define solutions to the nonlocal Cauchy problem. Note that since we allow for measure initial data and in particular, we want to include the fundamental solution in the solution concept. To this end, we need to work with a weaker notion of solutions compared to \cite[Definition 2.5]{LiWe24}, where the initial data were always assumed to be in $L^2(\R^d)$. In particular, we do not assume $\Vert u(t) \Vert_{L^2_{\loc}(\R^d)}$ to be bounded up to time $t=0$ nor $\Vert u(t) \Vert_{L^1(\R^d;w)}\in L^1(0,T)$.

\begin{definition}[Cauchy problem]
\label{def:Cauchy}
    Let $[\eta,T] \subset \R$ be an interval and $\mu \in \cM$ be a Radon measure. We say that 
    $u$ is a weak solution to the Cauchy problem
    \begin{align}
    \label{eq:Cauchy-measure}
        \begin{cases}
            \partial_t u - \mathcal{L}_t u &= 0 ~~ \text{ in } (\eta,T) \times \R^d,\\
            u(\eta) &= \mu ~~ \text{ in } \R^d,
        \end{cases}
    \end{align}
    if $u$ is a weak global solution in $(\eta+\varep, T)\times\R^d$ for any $\varep\in (0,T-\eta)$ in the sense of \autoref{def:global-sol}, and 
    \begin{align}\label{eq:initial-mu}
        \int_{\R^d} \psi(x) \d \mu(x) = \lim_{t \searrow \eta} \int_{\R^d} \psi(x) u(t,x) \d x, ~~ \forall\, \psi \in C_c(\R^d).
    \end{align}
\end{definition}
Recall that the dual space of $C_c(\R^d)$ is the space of all Radon measures and the convergence in \eqref{eq:initial-mu} is sometimes termed ``vague" convergence.

Finally, we introduce solutions to the dual Cauchy problem. 
\begin{definition}[dual Cauchy problem]
\label{def:dual-Cauchy}
Let $[\eta,T] \subset \R$ be an interval and $\mu \in \cM$ be a Radon measure. We say that $u$ is a weak solution to the dual Cauchy problem
    \begin{align}
    \label{eq:Cauchy-measure-dual}
        \begin{cases}
            \partial_t u + \mathcal{L}_t u &= 0 ~~ \text{ in } (\eta,T) \times \R^d,\\
            u(T) &= \mu ~~ \text{ in } \R^d,
        \end{cases}
    \end{align}
if $v(t,x):=u(T-t,x)$ is a solution to the Cauchy problem
    \begin{align*}
        \begin{cases}
            \partial_t v - \mathcal{L}_{T-t} v &= 0 ~~ \text{ in } (0,T-\eta) \times \R^d,\\
            v(0) &= \mu ~~ \text{ in } \R^d,
        \end{cases}
    \end{align*}
in the sense of \autoref{def:Cauchy}.
\end{definition}

\section{Weak Harnack inequalities}
\label{sec:wHI}

Let us collect some known results on the weak Harnack inequality for nonnegative supersolutions. All of these results have a local thrust.
\begin{theorem}[weak Harnack inequality]
\label{thm:weak-Harnack}
Assume that the kernel $K$ satisfies the upper bound \eqref{eq:up-bd} and the coercivity  condition \eqref{eq:coercive}.
Let $u \ge 0$ be a weak supersolution to \eqref{Eq:1:1} in $[t_o , t_o + 8R^{2s}] \times B_{8R}(x_o)$ in the sense of \autoref{def:local-sol}. Then, there exists $c>0$ depending on the data $\{d,s,\lm,\Lm\}$, such that
\begin{align}\label{eq:WH:1}
\bint_{t_o+R^{2s}}^{t_o+2R^{2s}}\bint_{B_R(x_o)} u(t,x) \d x \d t \le c \inf_{(t_o + 4R^{2s}, t_o + 8R^{2s}) \times B_R(x_o)} u.
\end{align}
Moreover, if the lower bound \eqref{eq:low-bd} of $K$ is satisfied, then 
\begin{align}\label{eq:WH:2}
\bint_{t_o+R^{2s}}^{t_o + 2R^{2s}} \bigg[ R^{2s} \int_{\R^d\setminus B_R(x_o)} \frac{u(t,x)}{|x-x_o|^{d+2s}} \d x \bigg] \d t \le c \inf_{(t_o + 4R^{2s}, t_o + 8R^{2s}) \times B_R(x_o)} u.
\end{align}
\end{theorem}

\begin{proof}
    The first claim \eqref{eq:WH:1} was shown in \cite{FeKa13}. The second claim  \eqref{eq:WH:2} was established in \cite[Theorem 1.9]{KaWe23b}.
\end{proof}

\begin{remark}
    Note that a weak Harnack inequality with a tail-term similar to \eqref{eq:WH:3} remains true if $K$ satisfies the upper bound \eqref{eq:up-bd}, the coercivity \eqref{eq:coercive}, and the UJS condition 
    \begin{align}
    \label{eq:UJS}
    K(t;x,y) \le \Lambda\, \bint_{B_r(x)} K(t;z,y) \d z
\end{align}
for  any $t\in\R$, $x,\,y\in\R^d$. and $r \le ( \frac{1}{4} \wedge \frac{|x-y|}{4} )$. In that case, if $u \ge 0$ is a weak solution in the sense of \autoref{def:local-sol} (instead of a weak supersolution)  we have the following integral estimate (see \cite[Proof of Theorem 6.4]{KaWe23b}):
\begin{align}\label{eq:WH:1.5}
\bint_{t_o+R^{2s}}^{t_o + 2R^{2s}} \bigg[ R^{2s} \sup_{y\in B_{\frac34R}(x_o)}\int_{\R^d\setminus B_R(x_o)} u(t,y) K(t;x,y) \d x \bigg] \d t \le c \inf_{(t_o + 4R^{2s}, t_o + 8R^{2s}) \times B_R(x_o)} u.
\end{align}
\end{remark}

The next result improves \eqref{eq:WH:1} and was shown in \cite[Theorem 1.5]{LiWe24}. It is of central importance in the theory of initial traces, which we establish in this article (see \autoref{lemma:ex-un-initial-trace}).

\begin{theorem}[improved weak Harnack inequality]
\label{thm:improved-weak-Harnack}
Assume that the kernel $K$ satisfies the upper bound \eqref{eq:up-bd} and the coercivity condition \eqref{eq:coercive}.
Let $u \ge 0$ be a weak supersolution to \eqref{Eq:1:1} in $[t_o , t_o + 8R^{2s}] \times B_{8R}(x_o)$ in the sense of \autoref{def:local-sol}. Then, there exists $c>0$ depending on the data $\{d,s,\lm,\Lm\}$, such that
\begin{align}\label{eq:WH:3}
\sup_{(t_o, t_o + R^{2s})} \bint_{B_R(x_o)} u(\cdot,x) \d x \le c \inf_{(t_o + 2R^{2s}, t_o + 8R^{2s}) \times B_R(x_o)} u.
\end{align}
\end{theorem}

\section{The fundamental solution and two-sided estimates}
\label{sec:fund-sol}

In the existing literature on a nonlocal Widder-type theory (see \cite{BPSV14,BSV17,Rui25,GQSV25}), the nonlocal kernels are always assumed to be translation invariant. This is in stark contrast to our article, where general nonlocal operators with bounded measurable coefficients are considered. A notable difference arising from our more general setup is that the fundamental solution cannot anymore be written explicitly or calculated via Fourier transform. In fact, its existence is already nontrivial.

In Subsection \ref{subsec:ex-fs} we complete the construction of the fundamental solution, which was initiated in \cite{LiWe24}. In particular, we prove that it is a solution to the Cauchy problem with a Dirac measure as the initial datum.

Moreover, in Subsection \ref{subsec-upper} and Subsection \ref{subsec-lower} we establish sharp two-sided estimates for the fundamental solution, thereby proving \autoref{thm:hkb}, a central result of this article.

\subsection{Existence of the fundamental solution}
\label{subsec:ex-fs}

For $\eta \in [0,T)$ consider the Cauchy problem
\begin{align}
\label{eq:inhom-Cauchy}
    \begin{cases}
        \partial_t u - \mathcal{L}_t u &= 0 ~~ \text{ in } (\eta,T) \times \R^d,\\
        u(\eta) &= f ~~ \text{ in } \R^d,
    \end{cases}
\end{align}
where $\mathcal{L}_t$ is a time-dependent operator satisfying \eqref{eq:up-bd} and \eqref{eq:coercive}, and $f \in L^2(\R^d)$. Analogously, for $\xi\in(0,T]$ we also consider the dual problem
\begin{align}
\label{eq:dual-inhom-Cauchy}
    \begin{cases}
        \partial_t u + \mathcal{L}_t u &= 0 ~~ \text{ in } (0,\xi) \times \R^d,\\
        u(\xi) &= f ~~ \text{ in } \R^d.
    \end{cases}
\end{align}
Solutions to \eqref{eq:inhom-Cauchy} and \eqref{eq:dual-inhom-Cauchy} are to be understood in the sense of \autoref{def:Cauchy} and \autoref{def:dual-Cauchy}, respectively, with $\mu = f$.

In \cite{LiWe24} we have  established a representation formula for solutions to the Cauchy problems \eqref{eq:inhom-Cauchy} and \eqref{eq:dual-inhom-Cauchy} with $L^2$ initial data. Let us summarize our results as follows.

\begin{proposition}
    \label{prop:representation}
    Assume that the kernel $K$ satisfies the upper bound \eqref{eq:up-bd} and the coercivity condition \eqref{eq:coercive}. Let $\eta \in [0,T)$. Then, there exists a function $(t,x,y) \mapsto p_{\eta,t}(x,y)$, such that for any $f \in L^2(\R^d)$ the function
    \begin{align}
    \label{eq:representation}
        (t,x) \mapsto P_{\eta,t}f(x) := \int_{\R^d} p_{\eta,t}(x,y) f(y) \d y
    \end{align}
    is the unique solution to \eqref{eq:inhom-Cauchy} in the sense of \autoref{def:Cauchy} which attains its initial datum in the sense of $L^2_{\loc}(\R^d)$. 
    Moreover, it holds for any $x,y \in \R^d$ and $0 \le \eta < \tau < t < T$ that 
    \begin{align}
    \label{eq:p-properties}
        p_{\eta,t}(x,y) \ge 0,  ~~ \int_{\R^d} p_{\eta,t}(x,y) \d y \le 1, ~~ p_{\eta,t}(x,y) = \int_{\R^d} p_{\tau,t}(x,z) p_{\eta,\tau}(z,y) \d z.
    \end{align}
   Likewise, for $0<t<\xi\le T$ there exists a function $(t,x,y) \mapsto \hat{p}_{\xi,t}(x,y)$, such that for any $f \in L^2(\R^d)$ the function
        \begin{align}
    \label{eq:dual-representation}
        (t,x) \mapsto \hat{P}_{\xi,t} f(x) := \int_{\R^d} \hat{p}_{\xi,t}(x,y) f(y) \d y
    \end{align}
    solves the dual problem \eqref{eq:dual-inhom-Cauchy}. Moreover,  the following symmetry holds:
    \begin{align}
    \label{eq:dual-relation}
        p_{\eta,t}(x,y) = \hat{p}_{t,\eta}(y,x) \qquad \forall\, x,y \in \R^d, ~~ \forall\, 0 \le \eta < t \le T.
    \end{align}
\end{proposition}

\begin{remark}
    Note that due to $f \in L^2(\R^d)$, weak solutions constructed in \autoref{prop:representation} possess stronger regularity than the one assumed in \autoref{def:Cauchy}. In fact, they belong to the space $L^{\infty}((\eta,T);L^{2}(\R^d)) \cap L^2((\eta,T);H^s(\R^d))$, see \cite[Definition~2.5]{LiWe24} in this connection.
\end{remark}

The density function $p_{\eta,t}(x,y)$ is the fundamental solution associated to the nonlocal parabolic operator $\partial_t-\mathcal{L}_t$, and in particular, it is a solution to the Cauchy problem with a Dirac measure as the initial datum. We will give a direct proof of this fact in the following lemma without relying on any semigroup theory. This is in stark contrast to the existing literature.

\begin{lemma}
\label{lemma:heat-kernel-sol-property}
    Assume that the kernel $K$ is as in \autoref{prop:representation}. Then, the function $(t,x) \mapsto p_{\eta,t}(x,y)$ solves the Cauchy problem \eqref{eq:Cauchy-measure} in $(\eta,T) \times \R^d$ in the sense of \autoref{def:Cauchy} with $\mu = \delta_{\{ y \}}$ for every $y \in \R^d$. Moreover, the function $(t,x) \mapsto \hat{p}_{t,\xi}(x,y)$ solves the dual Cauchy problem \eqref{eq:Cauchy-measure-dual} in $(0,\xi) \times \R^d$ in the sense of \autoref{def:dual-Cauchy} with $\mu = \delta_{\{ y \}}$ for every $y \in \R^d$.
\end{lemma}

\begin{proof}
    First, we claim that for any $\tau \in (\eta,T)$ and $y \in \R^d$ there holds
    \begin{align}
    \label{eq:p-in-L^2}
        f(x) := p_{\eta,\tau}(x,y) \in L^2(\R^d).
    \end{align}
    To see this, recall that by the proof of \cite[Lemma A.4]{LiWe24}, we have
    \begin{align*}
        \Vert \hat{P}_{\tau,\eta} g \Vert_{L^{\infty}(\R^d)} \le c(\eta,\tau) \Vert g \Vert_{L^2(\R^d)}, \qquad \forall\, g \in L^2(\R^d).
    \end{align*}
    In particular, by \eqref{eq:dual-relation}, we deduce
    \begin{align*}
        \left| \int_{\R^d} p_{\eta,\tau}(z,y) g(z) \d z \right| = |\hat{P}_{\tau,\eta} g(y)| \le c(\eta,\tau) \Vert g \Vert_{L^2(\R^d)},
    \end{align*}
    which implies \eqref{eq:p-in-L^2} by duality.
    
    By the semigroup property \eqref{eq:p-properties}, $(t,x) \mapsto p_{\eta,t}(x,y)$ is of the form \eqref{eq:representation}, i.e. it holds
    \begin{align*}
        p_{\eta,t}(x,y) = P_{\tau,t}\left( p_{\eta,\tau}(\cdot,y) \right)(x) = P_{\tau,t} f(x),
    \end{align*}
    and therefore we deduce from \autoref{prop:representation} that $(t,x) \mapsto p_{\eta,t}(x,y)$ is a global solution in $(\tau, T)\times\R^d$, in the sense of \autoref{def:global-sol}. By arbitrariness of $\tau\in(\eta, T)$, it is a global solution in $(\eta+\varepsilon, T)\times\R^d$, $\forall\,\varepsilon\in(0,T-\eta)$. 
    
    To see that $p_{\eta,t}(\cdot,y) \to \delta_{\{ y \}}$ as $t\searrow\eta$, we take $\psi \in C_c(\R^d)$ and compute
    \begin{align*}
        \lim_{t \searrow \eta} \int_{\R^d} \psi(x) p_{\eta,t}(x,y) \d x = \lim_{t \searrow \eta} (P_{\eta,t} \psi)(y) = \psi(y) = \int_{\R^d} \psi(x) \d \delta_{ \{y\} }(x),
    \end{align*}
    for every $y\in \R^d$, where we employed \autoref{prop:initial} in order to identify the limit in the second equality for every $y \in \R^d$. Note that $P_{\eta,t} \psi$ is globally bounded by $\|\psi\|_{L^{\infty}(\R^d)}$; see~\cite[Remark A.3]{LiWe24}. The results for $(t,x) \mapsto \hat{p}_{\xi,t}(x,y)$ hold true by analogous arguments.
\end{proof}

\begin{remark}
    As we shall see in \autoref{lemma:uniqueness},  the fundamental solution is uniquely characterized by the properties listed in \autoref{lemma:heat-kernel-sol-property}.
\end{remark}

\begin{remark}
    For the proof of the upper bound of the fundamental solution in \cite{KaWe23} it is also important to observe that \autoref{prop:representation} and \autoref{lemma:heat-kernel-sol-property} remain true for operators $\mathcal{L}_t$ with kernels $K(t;x,y) \chi_{\{|x-y| \le \rho\}}(x,y)$, where $\rho > 0$. Note that these kernels only satisfy the coercivity condition \eqref{eq:coercive} for $\rho \ge c r$, however, the statements of \autoref{prop:representation} and \autoref{lemma:heat-kernel-sol-property} remain true by the same proofs. In fact, the parabolic weak maximum principle which was used in the proof of \cite[Lemma A.2]{LiWe24} remains true (see \cite[Lemma 5.2]{KaWe23}) for the truncated kernels. Also the local boundedness which was used in the proof of \cite[Lemma A.4]{LiWe24} remains true for truncated kernels, at least for small enough balls (see also \cite[Lemma 2.4]{KaWe23b}), since the coercivity condition holds on small balls. The same argument allows one to reproduce the proof of the symmetry of the truncated fundamental solution by splitting into balls with small enough radii \cite[Section A.4]{LiWe24} and the proof of the continuity of solutions at the initial level (see \autoref{prop:initial}). 
\end{remark}

\subsection{Upper bounds of the fundamental solution}
\label{subsec-upper}

Heat kernel upper bounds were already established in \cite{KaWe23} by an entirely analytic technique, based on the original ideas by Aronson for second order operators. In \cite{KaWe23} the authors make several mild assumptions, on the existence of the fundamental solution, as well as on its symmetry and integrability properties. These assumptions were well-known to be true for time-independent kernels $K$, and their verification relied on semigroup theory. Using \autoref{prop:representation}, we can now verify these assumption also for time-dependent problems. Note that our proof only makes use of variational techniques. 

Hence, we have the following result:

\begin{proposition}\label{prop:up-bd}
    Assume that the kernel $K$ satisfies the upper bound \eqref{eq:up-bd} and the coercivity condition \eqref{eq:coercive}. Let $p_t(x,y)\equiv p_{0,t}(x,y)$ be identified in \autoref{prop:representation}.  There exists $c > 0$, depending only on the data $\{d,s,\lambda,\Lambda\}$, such that 
\begin{align}
\label{eq:heat-kernel-upper}
p_t(x,y) \le  c \left( t^{-\frac{d}{2s}} \wedge \frac{t}{|x-y|^{d+2s}}  \right)
\end{align}
holds for any $x,\,y\in\R^d$ and $t \in (0,T)$.
\end{proposition}

\begin{proof}
The proof of the on-diagonal upper bound, i.e.
\begin{align*}
    p_t(x,y) \le c\, t^{-\frac{d}{2s}}
\end{align*}
goes in the same way as in \cite[Theorem 2.3]{KaWe23}, using that the fundamental solution is a global solution in $(0,T)\times \R^d$, which we have verified in \autoref{lemma:heat-kernel-sol-property}.
The remaining necessary ingredients for the proof of the off-diagonal bound are listed in \cite[(2.6)--(2.9)]{KaWe23}, and they are all verified in \autoref{prop:representation}. The only difference is the symmetry assumption \cite[(2.7)]{KaWe23}, which needs to be replaced by \eqref{eq:dual-relation} in the proofs of \cite[Lemma 2.2, Theorem 3.4]{KaWe23}, and using that all results that are required in the proof also hold true for the dual equation.
\end{proof}

\subsection{Lower bounds of the fundamental solution}
\label{subsec-lower}

In this section, we aim to depict the lower bound of the fundamental solution.
In particular, we establish a lower bound that matches \eqref{eq:heat-kernel-upper} provided the lower bound \eqref{eq:low-bd} of the kernel is assumed. 

Such an estimate was previously only known for kernels that do not depend on $t$. Moreover, our proof is much shorter compared to the ones that are available in the literature and only uses the nonlocal weak Harnack inequality.

\begin{proposition}\label{prop:low-bd}
    Assume that the kernel $K$ satisfies the upper bound \eqref{eq:up-bd} and the lower bound \eqref{eq:low-bd}. Let $p_t(x,y)\equiv p_{0,t}(x,y)$ be identified in \autoref{prop:representation}.  There exists $c > 0$, depending only on the data $\{d,s,\lambda,\Lambda\}$, such that 
\begin{align*}
c\, p_t(x,y) \ge  t^{-\frac{d}{2s}} \wedge \frac{t}{|x-y|^{d+2s}} 
\end{align*}
holds for any $x,\,y\in\R^d$ and $t \in (0,T)$.
\end{proposition}

The proof relies on considering two situations.
First, we prove the following \textbf{near-diagonal} lower bound, which holds true under the weaker coercivity condition \eqref{eq:coercive} than \eqref{eq:low-bd}.

\begin{lemma}
\label{lemma:NDL}
Assume that the kernel $K$ satisfies the upper bound \eqref{eq:up-bd} and the coercivity condition \eqref{eq:coercive}. Let $p_t(x,y)\equiv p_{0,t}(x,y)$ be identified in \autoref{prop:representation}. There exists $c > 0$, depending only on the data $\{d,s,\lambda,\Lambda\}$, such that 
\begin{align*}
c\, p_t(x,y) \ge t^{-\frac{d}{2s}} 
\end{align*}
holds for any $x,\,y\in\R^d$ and $t>0$ satisfying $|x-y|^{2s}\le t$.
\end{lemma}

\begin{proof}
A proof was given in \cite[Theorem 19.1.1]{Wei22}, relying on the weak Harnack inequality \eqref{eq:WH:1}. Here, we include a slightly shorter argument based on the improved weak Harnack inequality \eqref{eq:WH:3}. Consider $R>0$ to be chosen and an arbitrary ball $B_{2R}(z)$. Applying the representation formula \eqref{eq:dual-representation} with $\xi=R^{2s}$ and $f=\chi_{B_{2R}(z)}$ we conclude that
\begin{align*}
        (t,x) \mapsto \int_{B_{2R}(z)} \hat{p}_{\xi,t}(x,y) \d y
\end{align*}
is the solution to the dual problem \eqref{eq:dual-inhom-Cauchy} in the sense of \autoref{def:dual-Cauchy}. After the change of variable $t\to \xi-t$, we may apply the weak Harnack inequality \eqref{eq:WH:3} to this solution and then exploit the symmetry property \eqref{eq:dual-relation} to obtain that there exists $c_o=c_o(d,s,\lambda,\Lambda)$, such that
\[
\int_{B_{2R}(z)} p_{0,R^{2s}}(y,x) \d y= \int_{B_{2R}(z)} \hat{p}_{R^{2s},0}(x,y) \d y\ge c_o
\]
holds for all $x\in B_{2R}(z)$. Since $(t,y)\mapsto p_{0,t}(y,x)$ is a solution to \eqref{eq:Cauchy-measure} by \autoref{lemma:heat-kernel-sol-property}, we may rely on the above integral estimate and apply \eqref{eq:WH:3} to $(t,y)\mapsto p_{0,t}(y,x)$ at $t_o=R^{2s}$ and conclude that for any $x\in B_{2R}(z)$, we have
\begin{align*}
\inf_{y\in B_{2R}(z)} p_{0,(2R)^{2s}}(y,x) \ge \frac{c}{R^d}
\end{align*}
for some $c=c(d,s,\lambda,\Lambda)$. Finally, we argue as follows. Given $t>0$ we let $(2R)^{2s}=t$ and require $x,\,y$ to satisfy $|x-y|\le 4R$. In this way, there exists some $B_{2R}(z)$ whose closure includes $x$ and $y$. Then, one realizes that the last display is exactly what has been desired after redefining $c$ in terms of $\{d,s,\lambda,\Lambda\}$.
\end{proof}

Next, we carry out a detailed study of the \textbf{off-diagonal} lower bound for the fundamental solution. Even though for the near-diagonal estimate, the coercivity \eqref{eq:coercive} suffices, it turns out the off-diagonal estimates need more particular information regarding the kernel than merely \eqref{eq:coercive}. The next result establishes a lower bound assuming the lower bound \eqref{eq:low-bd} of the kernel, and hence we match the upper bound \eqref{eq:heat-kernel-upper} of the fundamental solution as stated in \autoref{prop:low-bd}.

\begin{lemma}
\label{lemma:off-diag}
Assume that the kernel $K$ satisfies the upper bound \eqref{eq:up-bd} and the lower bound \eqref{eq:low-bd}. Let $p_t(x,y)\equiv p_{0,t}(x,y)$ be identified in \autoref{prop:representation}.  There exists $c > 0$, depending only on the data $\{d,s,\lambda,\Lambda\}$, such that 
\begin{align*}
c\, p_t(x,y) \ge  \frac{t}{|x-y|^{d+2s}} 
\end{align*}
holds for any $x,\,y\in\R^d$ and $t \in (0,T)$ satisfying  $|x-y|^{2s}\ge t$.
\end{lemma}

\begin{proof}
Fix $x,y \in \R^d$ and $t > 0$ with $|x-y|^{2s} \ge t$ and choose $R^{2s}=\frac14 t$. Then, using \autoref{lemma:heat-kernel-sol-property}, by \eqref{eq:WH:2} we have for some $c=c(d,s,\lambda,\Lambda)$ that
\begin{align*}
c\, p_t(x,y) &\ge 
 R^{2s} \bint_{t - 2R^{2s}}^{t - R^{2s}} \int_{\R^d \setminus B_R(y)} \frac{p_{\tau}(x,z)}{|y-z|^{d+2s}} \d z \d \tau \\
 &\ge R^{2s} \bint_{t - 2R^{2s}}^{t - R^{2s}} \int_{B_{R}(x)} \frac{p_{\tau}(x,z)}{|y-z|^{d+2s}} \d z \d \tau.
\end{align*}
Here, to obtain the last line we used $|x-y| \ge t^{\frac{1}{2s}}=4^{\frac{1}{2s}} R$  by construction, and hence $B_{R}(x) \subset \R^d \setminus B_R(y)$. Next, we aim to bound the integrand from below. To this end, consider $z \in B_{R}(x)$ and $\tau \in (t - 2R^{2s} , t - R^{2s})$ and estimate 
\begin{align*}
|x-z|^{2s} <  R^{2s} = \tfrac{1}{4}t \le t-2R^{2s}  \le  \tau .
\end{align*}
Therefore, by \autoref{lemma:NDL} there exists $c_o=c_o(d,s,\lambda,\Lambda)$, such that
\begin{align*}
p_{\tau}(x,z) \ge \frac{c_o}{\tau^{\frac{d}{2s}}} \ge \frac{c_o}{t^{\frac{d}{2s}}}.
\end{align*}
Moreover, for any $z \in B_{R}(x)$ it holds
\begin{align*}
|y-z| \le |x-y| + |x - z| \le |x-y| + R \le 2|x-y| .
\end{align*}
Thus, we estimate
\begin{align*}
c\, p_t(x,y) &\ge  R^{2s} \bint_{t - 2R^{2s}}^{t - R^{2s}} \int_{B_{R}(x)} \frac{p_{\tau}(x,z)}{|y-z|^{d+2s}} \d z \d \tau \\
&\ge \frac{c_o}{t^{\frac{d}{2s}}} \frac{ R^{2s}}{|x-y|^{d+2s}} \bint_{t - 2R^{2s}}^{t - R^{2s}} \int_{B_{R}(x)}  \d z \d \tau \\
&\ge  c_o \frac{R^{2s}}{|x-y|^{d+2s}} \frac{R^{d}}{t^{\frac{d}{2s}}}  \\
&\ge c_o \frac{t}{|x-y|^{d+2s}}
\end{align*}
where $c$ and $c_o$ are generic constants depending on $\{d,s,\lambda,\Lambda\}$.
This concludes the proof upon redefining $c/c_o$ as $c$.
\end{proof}

\begin{proof}[Proof of \autoref{prop:low-bd}]
    Combining \autoref{lemma:NDL} and \autoref{lemma:off-diag} yields the desired result.
\end{proof}

\begin{remark}
    \label{rem:gen-lb}
    It is an intriguing question to investigate the validity of lower heat kernel estimates for kernels that violate the pointwise lower bound \eqref{eq:low-bd}, and so far, sharp lower bounds are only available in special cases. For instance in \cite{Kan23}, relying on probabilistic arguments, the author investigates kernels that are comparable to $|x-y|^{-d-2s}$ only on a fixed double cone centered at $x-y$ and otherwise zero.
    In such a case, \eqref{eq:low-bd} is violated while the significantly weaker UJS condition \eqref{eq:UJS} is in force. 
    
    A straightforward adaptation of the proof of \autoref{lemma:off-diag}, but using \eqref{eq:WH:1.5} instead of \eqref{eq:WH:2}, also allows us to prove the following off-diagonal lower bound for the fundamental solution in case $K$ satisfies only \eqref{eq:up-bd}, \eqref{eq:coercive}, and \eqref{eq:UJS}, whenever $|x-y|^{2s} \ge t$:
    \begin{align*}
        p_t(x,y) &\ge \frac{c_o}{t^{\frac{d}{2s}}} R^{2s}  \bint_{t - 2R^{2s}}^{t - R^{2s}} \sup_{w \in B_{\frac{3}{4}R}(y)} \int_{B_{R}(x)} K(\tau;z,w) \d z \d \tau \\
        &\ge \frac{c_o}{t^{\frac{d}{2s}}} R^{d+2s}  \bint_{t - 2R^{2s}}^{t - R^{2s}} K(\tau;x,y)  \d \tau \\
        &\ge c_o\, t  \bint_{\frac{1}{2}t}^{\frac{3}{4}t} K(\tau;x,y)  \d \tau.
    \end{align*}
    Hence, we obtain for any $x,y \in \R^d$ that
    \begin{align}
    \label{eq:lower-UJS-heatkernel}
        c\, p_t(x,y) \ge t^{-\frac{d}{2s}} \wedge t  \bint_{\frac{1}{2}t}^{\frac{3}{4}t} K(\tau;x,y)  \d \tau.
    \end{align}
    Note that a similar estimate has been established in \cite[Theorem 4.8]{GHH17} via a completely different approach restricted to time-independent kernels. Moreover, we emphasize that \eqref{eq:lower-UJS-heatkernel} is sufficient to deduce the lower heat kernel bound in \cite{Kan23} by decomposing the space in the same way.
\end{remark}

\section{Widder-type theorem}
\label{sec:Widder}

The goal of this section is to prove the Widder-type theorem for solutions to the nonlocal Cauchy problem (see \autoref{thm:widder}).

The proof relies on two main steps. First, we establish that any nonnegative solution in $(0,T) \times \R^d$ is a solution to the Cauchy problem with some unique initial datum $\mu \in \cM_+$ satisfying the growth condition \eqref{eq:mu-growth} (see Subsection \ref{subsec:initial-trace}). In a second step, we show that given any initial datum $\mu \in \cM_+$ with such growth condition, there exists a unique solution to the Cauchy problem (see Subsection \ref{sec:ex-un}).

\subsection{A weighted estimate}
Prior to identifying the initial trace of global solutions, we need a technical lemma based on the improved weak Harnack inequality (see \autoref{thm:improved-weak-Harnack}).
\begin{lemma}
\label{lemma:weighted-L1}
    Assume that the kernel $K$ satisfies the upper bound \eqref{eq:up-bd} and the coercivity condition \eqref{eq:coercive}.
    Let $u$ be a nonnegative, global weak solution in $(\varep,T]\times \R^d$, $\forall\,\varep\in(0,T)$ in the sense of \autoref{def:global-sol}.  Then, there exists $c > 0$, depending only on the data $\{d, s, \lambda, \Lambda\}$ and $T$, such that for any $0 < t < T$, 
    \begin{align*}
        \Vert u(t) \Vert_{L^1(\R^d;w)} \le c\, \Vert u(T) \Vert_{L^1(\R^d;w)}.
    \end{align*}
\end{lemma}
\begin{proof}
Suppose first that $T\le 1$ and that $\Vert u(T) \Vert_{L^1(\R^d;w)}<\infty$ without loss of generality.
Let us denote $r=(T-t)^{\frac1{2s}}$ and decompose $\R^d$ into pairwise disjoint subcubes $K_r(x_o)$ with side length $2r$, center $x_o$ and faces parallel with the coordinate planes.  Apply the improved weak Harnack inequality \eqref{eq:WH:3} to each of them and obtain that
\[
\int_{K_r(x_o)} u(t,x) \d x \le c \int_{K_r(x_o)} u(T,x)\d x
\]
for some $c=c(d,s,\lm, \Lm)$. 
Since for any $x \in K_r(x_o)$ with $r<1$, it holds that
$$
\frac{2^{-(d+2)}}{(1+|x|)^{d+2s}} \le 
\frac1{(1+|x_o|)^{d+2s}} \le \frac{2^{d+2}}{(1+|x|)^{d+2s}},
$$ 
 we also have a weighted version of the previous integral estimate:
\begin{align*}
    \int_{K_r(x_o)} \frac{u(t,x)}{(1+|x|)^{d+2s}} \d x \le c \int_{K_r(x_o)} \frac{u(T,x)}{(1+|x|)^{d+2s}} \d x
\end{align*}
for some $c=c(d,s,\lm,\Lm)$.
Summing  over all subcubes that constitute $\R^d$ on both sides gives the desired result, under the assumption that $T\le 1$. The general case is proven by slicing $(0,T)$ into subintervals and repeating the same argument.
\end{proof}

\begin{remark}
    Note that a generalized version of the previous lemma is given by the time-insensitive $L^{1}(\R^d;w)$ estimate from \cite[Proposition~6.1]{LiWe24}, whose proof is significantly more involved. 
    In contrast, note that the proof of the previous lemma is merely based on local properties of supersolutions, and that the result is sufficient for our purpose.
\end{remark}

\subsection{Existence and uniqueness of initial trace}
\label{subsec:initial-trace}
Every global solution admits a Radon measure as its initial trace. Moreover, we can quantify its growth at large.

\begin{proposition}
\label{lemma:ex-un-initial-trace}
    Assume that the kernel $K$ satisfies the upper bound \eqref{eq:up-bd} and the coercivity condition \eqref{eq:coercive}. Let $u$ be a nonnegative, global weak solution in $(\varep,T] \times \R^d$, $\varep\in(0,T)$ in the sense of \autoref{def:global-sol}. Assume that $u(T) \in L^1(\R^d;w)$. Then, there exists a unique $\mu \in \cM_+$ with 
    \begin{align*}
        \int_{\R^d} \frac{\d \mu(x)}{(1 + |x|)^{d+2s}} \le c\, \Vert u(T) \Vert_{L^1(\R^d;w)},
    \end{align*}
    where $c > 0$ depends only on the data $\{d,s,\lm,\Lm\}$ and $T$, such that $u(0) = \mu$ in the sense of \eqref{eq:initial-mu}.
\end{proposition}

\begin{proof}
    By \autoref{lemma:weighted-L1}, there exists $c > 0$, depending only on $\{d, s, \lambda, \Lambda\}$ and $T$, such that
    \[
    \int_{B_R}u(t,x)\d x\le c\,(1+R)^{d+2s}\Vert u(T) \Vert_{L^1(\R^d;w)}
    \]
    holds true for any $B_R\equiv B_R(0)$ with $R>0$ and any $t\in(0,T)$. Therefore, $\{u(t,x)\d x: t\in(0,T)\}$ is a family of Radon measures on $\R^d$ that are locally, uniformly bounded. By compactness for measures with respect to vague convergence in the sense of \eqref{eq:initial-mu} (see \cite[Theorem 1.41]{EvGa15}), there exist a Radon measure $\mu$ and a sequence $s_n\searrow0$, such that
    \[
    \lim_{n\to\infty}\int_{\R^d} \vp (x) u(s_n,x)\d x=\int_{\R^d}\vp(x)\d\mu(x)\qquad\forall\,\vp\in C_c(\R^d).
    \]
    Moreover, it is not hard to see from the above uniform bound that
    \[
    \mu(B_R)\le c\,(1+R)^{d+2s}\Vert u(T) \Vert_{L^1(\R^d;w)},
    \]
    and hence the claimed estimate on the growth of $\mu$ follows.

    Next, we need to verify that the identification of $\mu$ is independent of the selection of sequences in order to deduce that 
    \[
    \lim_{t\searrow0}\int_{\R^d} \vp (x) u(t,x)\d x=\int_{\R^d}\vp(x)\d\mu(x)\qquad\forall\,\vp\in C_c(\R^d).
    \]
    To this end, we first prove that for any $t_1,\,t_2\in(0,\frac14 T)$, $R>0$ and $\sig\in(0,1)$, there exists $c>0$ depending on $\{d,s,\lm,\Lm\}$ and $T$, such that
    \begin{align}\label{eq:t1-t2}
        \int_{B_R} u(t_1,x) \d x - \int_{B_{(1+\sig)R}} u(t_2,x) \d x\le \frac{c\, R^d}{\sig}\Vert u(T) \Vert_{L^1(\R^d;w)}\, |t_1-t_2|.
    \end{align}
    
    It suffices to prove it in case $t_1<t_2$. We introduce a sequence of functions $\zeta_k\in C_c^{1}(0,\frac14 T)$, $k\in\N$ such that $\zeta_k \to \chi_{[t_1, t_2]}$ and $\partial_t \zeta_k \to \delta_{t_1} - \delta_{t_2}$. For $B_R\equiv B_R(0)$ and $\sig\in(0,1)$, let $\xi \in C_c^{1}(B_{(1+\sig)R})$ be such that $0 \le \xi \le 1$, $\xi \equiv 1$ in $B_R$, and $|\nabla \xi| \le 4 (\sig R)^{-1}$. Using $\z_k \xi$ as a testing function in the weak formulation of $u$ and letting $k\to\infty$, we obtain that
    \begin{align*}
        \int_{\R^d} \xi u(t,x) \d x \bigg|_{t_1}^{t_2} \le \int^{t_2}_{t_1}\int_{\R^d}\int_{\R^d} |u(t,x)-u(t,y)||\xi(x)-\xi(y)|K(t;x,y)\d x\d y\d t.
    \end{align*}
    Evoking  the properties of $\xi$, the right-hand side is estimated by
    \begin{align*}
    \mathbf{I}_R :=  4 \int_{t_1}^{t_2} \iint_{(B_{2R}^c\times B_{2R}^c)^c} |u(t,x) - u(t,y)|    (1 \wedge (\sig R)^{-1}|x-y|) K(t;x,y) \d y \d x \d t.
    \end{align*}
    We split the domain of integration as $ (B_{4R}\times B_{4R}) \cup (B_{2R}\times B^c_{4R})\cup (B^c_{4R} \times B_{2R}) \supset (B_{2R}^c\times B_{2R}^c)^c$ and then use the upper bound of $K$ and the symmetry of the integrals.
This results in
\begin{align*}
    \mathbf{I}_{R} &\le \int_{t_1}^{t_2} \int_{B_{4R}}   \int_{B_{4R}} |u(t,x) - u(t,y)| (1 \wedge (\sig R)^{-1}|x-y|) K(t;x,y)  \d y \d x \d t \\
    & \quad + 2\int_{t_1}^{t_2} \int_{B_{2R}} \int_{B_{4R}^c}  \frac{u(t,x)}{|x-y|^{d+2s}}  \d y \d x \d t \\
    &\quad + 2\int_{t_1}^{t_2} \int_{B_{4R}^c}   \int_{B_{2R}}  \frac{u(t,x)}{|x-y|^{d+2s}}  \d y \d x \d t \\
    &=: \textbf{I}_1 + \mathbf{I}_{2} + \mathbf{I}_{3}.
\end{align*}
To estimate $\mathbf{I}_{1}$, we employ \cite[Propositions~4.1, 4.2]{LiWe24} up to a proper scaling and obtain 
\begin{align*}
    \mathbf{I}_{1} &\le \frac{c\, R^d}{\sig} \int_{t_1}^{t_2+3(t_2-t_1)}\int_{B_{16R}}\frac{u(t,x)}{R^{d+2s}}\d x \d t\\
    &\le \frac{c\, R^d}{\sig} \int_{t_1}^{t_2+3(t_2-t_1)} \|u(t)\|_{L^1(\R^d;w)}\d t\\
    &\le \frac{c\, R^d}{\sig} \|u(T)\|_{L^1(\R^d;w)} |t_1-t_2|
\end{align*}
for some $c=c(d,s,\lm,\Lm,T)$. Here, in the last line we employed \autoref{lemma:weighted-L1}.
For $\mathbf{I}_{2}$, we simply estimate
\begin{align*}
        \mathbf{I}_{2} &\le \int_{t_1}^{t_2} \int_{B_{2R}} u(t,x) \bigg(\int_{B_{4R}^c} \frac{\d y}{|x-y|^{d+2s}}  \bigg) \d x \d t \\
        &\le c\, R^d \int_{t_1}^{t_2} \int_{B_{2R}} \frac{u(t,x)}{R^{d+2s}}\d x \d t\\
        &\le c\, R^d \|u(T)\|_{L^1(\R^d;w)}|t_1-t_2|
    \end{align*}
for some $c=c(d,s,\lm,\Lm,T)$.
 For $\mathbf{I}_{3}$, we use that when $y \in B_{2R}$ and $x \in B_{4R}^c$, then $|x| \le |x-y| + |y| \le |x-y| + 2R \le 2 |x-y|$ and thus
    \begin{align*}
        \mathbf{I}_{3} &\le c \, R^d\int_{t_1}^{t_2} \int_{B_{4R}^c}\frac{u(t,x)}{|x|^{d+2s}} \d y \d t \\
        &\le c\, R^d \|u(T)\|_{L^1(\R^d;w)}|t_1-t_2|
    \end{align*}
for some $c=c(d,s,\lm,\Lm,T)$. Collecting these estimates we have proven \eqref{eq:t1-t2}.

Now, suppose there are another sequence $\tau_m\searrow0$ and another Radon measure $\nu$, such that
\[
    \lim_{m\to\infty}\int_{\R^d} \vp (x) u(\tau_m,x)\d x=\int_{\R^d}\vp(x)\d\nu(x)\qquad\forall\,\vp\in C_c(\R^d).
\]
In \eqref{eq:t1-t2}, we substitute $(t_1,t_2)$ by the sequences $(s_n,\tau_m)$; we first let $n\to\infty$ and then $m\to\infty$; by \cite[Theorem~1.40~(iii)]{EvGa15}  we obtain
\[
\mu(B_R)\le \nu(B_{(1+\sig)R}).
\]
Using the arbitrariness of $\sig\in(0,1)$ and interchanging the role of $\mu$ and $\nu$ yields the equality $\mu(B_R)= \nu(B_{R})$ for any $R>0$. The same equality holds for any ball $B_R(x_o)$. To see this, we only need to apply a change of variable $x\to x-x_o$ in \eqref{eq:t1-t2} and realize that
\[
\int_{\R^d}\frac{u(T,x)}{(1+|x-x_o|)^{d+2s}}\d x<\infty.
\]
Then, the same argument as above gives the conclusion. As a result, we have $\mu=\nu$. 
\end{proof}

\subsection{Uniqueness of solutions for given nonnegative initial data}
\label{sec:ex-un}
Now we prove the uniqueness result among nonnegative solutions.
\begin{proposition}
\label{lemma:uniqueness}
    Assume that the kernel $K$ satisfies the upper bound \eqref{eq:up-bd} and the coercivity condition \eqref{eq:coercive}. Let $u_1$ and $u_2$ be nonnegative, weak solutions to the Cauchy problem in $(0,T] \times \R^d$ in the sense of \autoref{def:Cauchy}, with $u_1(0) = \mu = u_2(0)$ for some $\mu \in \cM_+$. Then, we have $u_1 = u_2$ in $(0,T] \times \R^d$.
\end{proposition}

\begin{proof}
    We define $u = u_1 - u_2$ and observe that $u$ solves the Cauchy problem in $(0,T] \times \R^d$ with $u(0) = 0$ in the sense of \eqref{eq:initial-mu}. It suffices to show that $u(T,\cdot) \equiv 0$. Without loss of generality we assume $\Vert u(T) \Vert_{L^1(\R^d;w)}<\infty$.
    Let us take a function $\psi \in C_c^{\infty}(\R^d)$ and define $\phi$ to be the weak solution to the dual Cauchy problem in $(0,T)$ with $\phi(T,\cdot) = \psi$ in the sense of \autoref{def:dual-Cauchy}. Suppose $\supp(\psi) \subset [0,T] \times B_R$ for some positive $R$. Due to the representation formula~\eqref{eq:dual-representation} it is not hard to see that for any $(t,x) \in [0,T] \times \R^d$:
    \begin{align}
    \label{eq:phi-upper}
        |\phi(t,x)| \le \Vert \psi \Vert_{L^{\infty}(\R^d)} \int_{B_R} \hat{p}_{T,t}(x,y) \d y \le \frac{c}{(1 + |x|)^{d+2s}},
    \end{align}
    where $c > 0$ depends on $d$, $s$, $\lambda$, $\Lambda$, $R$, $T$, and $\Vert \psi \Vert_{L^{\infty}(\R^d)}$.
    Indeed, after a simple change of variables $t\to T-t$, we may use \eqref{eq:heat-kernel-upper} to bound $\hat{p}_{T,t}(x,y)$. Hence, for any $x \in \R^d$ we first estimate 
    \begin{align*}
        |\phi(t,x)| &\le c\, \Vert \psi \Vert_{L^{\infty}(\R^d)} \bigg( \int_{ B_{(T-t)^{1/2s}}(x)} \frac{\d y}{(T-t)^{\frac{d}{2s}}}  + \int_{\R^d \setminus B_{(T-t)^{1/2s}}(x)} \frac{T-t}{|x-y|^{d+2s}} \d y \bigg)\\
        &\le c\, \Vert \psi \Vert_{L^{\infty}(\R^d)},
    \end{align*}
    where $c$ only depends on $d$, $s$, $\lambda$, and $\Lambda$.
    On the other hand, when $|x| > 2R$, we observe that $|x-y| \ge \frac12 |x|$ for any $y \in B_R$, and use again \eqref{eq:heat-kernel-upper} to further estimate
    \begin{align*}
        |\phi(t,x)| &\le c\,\Vert \psi \Vert_{L^{\infty}(\R^d)} \int_{B_R} \frac{T-t}{|x|^{d+2s}} \d y \\
        &\le c\,\Vert \psi \Vert_{L^{\infty}(\R^d)}  \frac{T|B_R|}{|x|^{d+2s}}. 
    \end{align*}
    This proves \eqref{eq:phi-upper}.

The following argument is modeled on the proof of \cite[Proposition~6.1]{LiWe24}. Let $\xi_R \in C_c^{1}(B_{2R})$ be such that $0 \le \xi_R \le 1$, $\xi_R \equiv 1$ in $B_R$, and $|\nabla \xi_R| \le 4 R^{-1}$. Moreover,  consider arbitrary $\tau_o,\,\tau_1$ satisfying $0< \tau_o < \tau_1 < T$. We introduce a sequence of functions $\zeta_k\in C_c^{1}(0,T)$, $k\in\N$ such that
\begin{equation*}
\left\{
    \begin{array}{cc}
        \zeta_k \to \chi_{[\tau_o, \tau_1]},  \quad
         \partial_t \zeta_k \to \delta_{\tau_o} - \delta_{\tau_1}\quad \text{as}\>\> k\to\infty , \\ [5pt]
         \displaystyle\supp(\zeta_k) \subset \Big[\tau_o, \frac{T +\tau_1}{2}\Big].
    \end{array}\right.
\end{equation*}

Modulo a proper time mollification, we observe 
\begin{align*}
        -\int_0^{T} \partial_t \zeta_k(t) \int_{\R^d} \xi_R(x) u(t,x) \phi(t,x) \d x \d t  &=   \int_0^{T} \zeta_k(t) \int_{\R^d} \xi_R(x) \partial_t u(t,x) \phi(t,x) \d x \d t\\
    &\quad + \int_0^{T} \zeta_k(t) \int_{\R^d} \xi_R(x) \partial_t \phi(t,x) u(t,x)  \d x \d t .
\end{align*}
To estimate the right-hand side, we test the weak formulation of $u$ with $\zeta_k\xi_R \phi$ and the weak formulation of $\phi$ with $\zeta_k\xi_R u$, and we obtain
\begin{align*}
    \int_0^{T} \zeta_k(t) \int_{\R^d} \xi_R(x) \partial_t u(t,x) \phi(t,x) \d x \d t&=-\int_0^{T} \zeta_k(t)   \cE^{(t)}(u(t) , \xi_R\phi(t)) \d t,\\
    \int_0^{T} \zeta_k(t) \int_{\R^d} \xi_R(x) \partial_t \phi(t,x) u(t,x)  \d x \d t&=\int_0^{T} \zeta_k(t)   \cE^{(t)}(\phi(t) , \xi_R u(t) )  \d t.
\end{align*}
To continue we plug these two identities back into the previous equation, recall the definition of $\cE^{(t)}$, and employ the following algebraic identity
    \begin{align*}
    \nonumber
        (\phi_1 - \phi_2)(u_1 \xi_1 - u_2 \xi_2) &= \frac{1}{2}(\phi_1 - \phi_2) (u_1 - u_2)(\xi_1 + \xi_2) + \frac{1}{2}(\phi_1 - \phi_2) (u_1 + u_2)(\xi_1 - \xi_2) \\ \nonumber
        &= (u_1 - u_2)(\phi_1 \xi_1 - \phi_2 \xi_2) - \frac{1}{2}(u_1 - u_2)(\xi_1 - \xi_2)(\phi_1 + \phi_2) \\ 
        &\quad+ \frac{1}{2}(\phi_1 - \phi_2)(u_1 + u_2)(\xi_1 - \xi_2).
    \end{align*}
As a result, we arrive at
\begin{align*}
     \bigg| \int_0^{T} &\partial_t \zeta_k(t) \int_{\R^d} \xi_R(x) u(t,x) \phi(t,x) \d x \d t \bigg| \\
    & \le \frac{1}{2}\bigg| \int_0^{T} \zeta_k(t) \bigg( \int_{\R^d} \int_{\R^d} (u(t,x) - u(t,y))(\phi(t,x) + \phi(t,y)) \\
    &\qquad\qquad\qquad\qquad\qquad\qquad\cdot (\xi_R(x) - \xi_R(y)) K(t;x,y) \d y \d x \bigg) \d t \bigg| \\
    &\quad \quad + \frac{1}{2} \bigg| \int_0^{T} \zeta_k(t) \bigg( \int_{\R^d} \int_{\R^d} (u(t,x) + u(t,y))(\phi(t,x) - \phi(t,y))\\
    &\qquad\qquad\qquad\qquad\qquad\qquad\cdot(\xi_R(x) - \xi_R(y)) K(t;x,y) \d y \d x \bigg)\d t \bigg| \\
    &=: \textbf{I}_R + \textbf{II}_R.
\end{align*}

 Letting first $k\to\infty$ and then $R\to\infty$, the left-hand side tends to
\begin{align*}
    \bigg| \int_{\R^d} u(\tau_o , x)  \phi(\tau_o,x) \d x - \int_{\R^d} u(\tau_1 , x)  \phi(\tau_1,x) \d x  \bigg|,
\end{align*}
whereas the right-hand side tends to $0$. Indeed, given \eqref{eq:phi-upper} these two terms will converge to $0$ as $R\to\infty$ uniformly in $k$. The argument is similar to those in \cite[Proposition~6.1]{LiWe24} upon further replacing $u$ in $\textbf{I}_R$ and $ \textbf{II}_R$ by $u_1-u_2$. To be precise, for $\textbf{I}_R$ we use $|u(t,x)-u(t,y)|\le |u_1(t,x)-u_1(t,y)|+|u_2(t,x)-u_2(t,y)|$ to split the integral into two and then apply \cite[Lemma 6.3]{LiWe24} with $(0,T)$ replaced by $(\tau_o,T)$; for $\textbf{II}_R$ we use $|u(t,x)-u(t,y)|\le u_1(t,x)+u_1(t,y)+u_2(t,x)+u_2(t,y)$ to split the integral into four and then apply \cite[Lemma 6.4]{LiWe24} with $(t_o,T)$ replaced by $(\tau_o,T)$. Therefore, we have
\begin{align*}
     \int_{\R^d} u(\tau_o , x)  \phi(\tau_o,x) \d x = \int_{\R^d} u(\tau_1 , x)  \phi(\tau_1,x) \d x  .
\end{align*}
The left-hand side in the above display vanishes as we let $\tau_o\searrow 0$, since the initial condition of $u$ ensures
\begin{equation}\label{eq:tau-to-0}
    \lim_{t\searrow0}\underbrace{\int_{\R^d} u(t , x) \phi(t,x) \d x}_{=:\mathbf{III}}=0.
\end{equation}
Assuming this for the moment, we then let $\tau_1\nearrow T$ and recall $\phi(T)=\psi$.  Consequently, we immediately deduce from the dominated convergence theorem and $u (T,\cdot) \in L^1(\R^d;w)$ that
    \begin{align*}
        \int_{\R^d} \psi(x) u(T,x) \d x = 0
    \end{align*}
    for any  $\psi \in C_c^{\infty}(\R^d)$. As a result of density, we must have $u(T) \equiv 0$. 
    
    Therefore, it remains to verify \eqref{eq:tau-to-0}. Let us rewrite it as
    \begin{align*}
        \mathbf{III} &=    \int_{\R^d}  u(t,x) \phi (0,x) \d x +  \int_{\R^d}  u(t,x) (\phi(t,x) - \phi (0,x)) \d x \\
        &=: \mathbf{III}_{1} + \mathbf{III}_{2}.
    \end{align*}
    We first deal with $\mathbf{III}_{1}$. In fact, since $u(0) = 0$, we can take $R=\sigma=1$, send $t_2\searrow0$, and relabel $t_1$ as $t$ in \eqref{eq:t1-t2} to obtain that for any $t\in(0,\frac14T)$,
    \begin{align*}
        \int_{K_1} u(t,x) \d x \le c\,t\Vert u(T) \Vert_{L^1(\R^d;w)}.
    \end{align*}
Here, we used that $\int_{K_2} u(t_2,x) \d x \to 0$ by \cite[Theorem 1.40 (iii)]{EvGa15}, and, with no loss of generality, we replaced the ball $B_1$ by the cube $K_1$. Up to a translation, we get from the previous estimate that 
    \begin{align*}
        \int_{K_1(x_o)} u(t,x) \d x \le c\,t\int_{\R^d} \frac{u(T,x)}{(1+|x-x_o|)^{d+2s}}\d x 
    \end{align*}
    for any $x_o\in\R^d$.
    Now that for any $x \in K_1(x_o)$, it holds that
$$
\frac{2^{-(d+2)}}{(1+|x|)^{d+2s}} \le 
\frac1{(1+|x_o|)^{d+2s}} \le \frac{2^{d+2}}{(1+|x|)^{d+2s}},
$$ 
 we also have a weighted version of the previous integral estimate:
 \begin{align*}
        \int_{K_1(x_o)} \frac{u(t,x)}{(1+|x|)^{d+2s}} \d x \le \frac{c\,t}{(1+|x_o|)^{d+2s}}\int_{\R^d} \frac{u(T,x)}{(1+|x-x_o|)^{d+2s}}\d x 
    \end{align*}
 for some $c>0$ depending on $\{d,s,\lm,\Lm\}$ and $T$. Summing over all cubes that constitute $\R^d$, we end up with 
 \begin{align*}
        \int_{\R^d} \frac{u(t,x)}{(1+|x|)^{d+2s}} \d x &\le c\, t\int_{\R^d}\frac{1}{(1+|y|)^{d+2s}}\int_{\R^d} \frac{u(T,x)}{(1+|x-y|)^{d+2s}}\d x  \d y\\
        &\le c\, t \int_{\R^d} \frac{u(T,x)}{(1+|x|)^{d+2s}}\d x.
    \end{align*}
Here, we used Fubini's theorem and the fact that
\begin{align*}
    \int_{\R^d}\frac{1}{(1+|y|)^{d+2s}}\frac{\d y}{(1+|x-y|)^{d+2s}} & \le \int_{B_1^c} \frac{1}{(1+|y|)^{d+2s}}\frac{\d y}{(1+|x-y|)^{d+2s}}  + \frac{c}{(1+|x|)^{d+2s}} \\
 & \quad\le  \frac{c}{(1+|x|)^{d+2s}}
\end{align*}
for some $c=c(d,s)$, which follows from \cite[Lemma~3.4]{LiWe24} and the observation that the estimate for the integral over $B_1$ is trivial.
    As a result of the above estimate, $\mathbf{III}_{1}$ converges to $0$ recalling \eqref{eq:phi-upper}. 
    For $\mathbf{III}_{2}$, we first claim that for some $c>0$ and $\alpha \in (0,1)$ depending only on $\{d,s,\lm,\Lm\}$, such that for any $t\in(0,1)$ and $x\in\R^d$ we have
    \begin{align}
    \label{eq:Holder-est-phi}
        |\phi(t,x) - \phi (0,x)| \le \frac{c\, t^{\alpha}}{(1 + |x|)^{d+2s}}.
    \end{align}
    This is a consequence of the H\"older regularity estimate (see \cite[Theorem~1.5]{KaWe23b}) for $\phi$. In fact, we estimate
    \begin{align*}
        |\phi(t,x) - \phi (0,x)| 
        \le c\, t^{\alpha} \bigg(\Vert \phi \Vert_{L^{\infty}((0,T) \times B_1(x))} +  \sup_{0<t<T}\int_{\R^d\setminus B_1(x)}\frac{\phi(t,y)}{|y-x|^{d+2s}}\d y\bigg).
    \end{align*}
    Clearly by \eqref{eq:phi-upper} the first term in the parenthesis is bounded by
    \begin{align*}
        \Vert \phi \Vert_{L^{\infty}((0,T) \times B_1(x))}\le c \,\Vert (1 + |\cdot|)^{-d-2s} \Vert_{L^{\infty}(B_1(x))} \le c\, (1 + |x|)^{-d-2s},
    \end{align*}
    whereas the second term admits a similar bound due to \cite[Lemma~3.4]{LiWe24}. Therefore, we can use \eqref{eq:Holder-est-phi} together with \autoref{lemma:weighted-L1} to estimate
    \begin{align*}
    \mathbf{III}_{2}&\le c\, t^{\al}\int_{\R^d}\frac{ u(t,x)}{(1+|x|)^{d+2s}}\d x\\
    &\le c\,t^{\al}\|u(T,\cdot)\|_{L^1(\R^d;w)}\to0
    \end{align*}
    as $t\searrow0$. This establishes \eqref{eq:tau-to-0} and completes the proof.
 \end{proof}

\subsection{Existence of solutions for given intitial data}
Next, we show that if a Radon measure satisfies a proper growth condition, then its ``convolution" with the fundamental solution gives a solution to the Cauchy problem with the same Radon measure as initial datum.
\begin{proposition}
\label{lemma:existence}
    Assume that the kernel $K$ satisfies the upper bound \eqref{eq:up-bd} and the coercivity condition \eqref{eq:coercive}. Let $\mu \in \cM$ satisfy  
    $$
    \int_{\R^d} \frac{\d |\mu|(x)}{(1 + |x|)^{d+2s}}  < \infty
    $$ 
    and let $p_t(x,y)\equiv p_{0,t}(x,y)$ be identified in \autoref{prop:representation}.
    Then, the function given by
    \begin{align*}
        u(t,x) := \int_{\R^d} p_t(x,y) \d \mu(y)
    \end{align*}
    solves the Cauchy problem in $(0,T) \times \R^d$ in the sense of \autoref{def:Cauchy}, with $u(0) = \mu$.
\end{proposition}

\begin{proof}
    Clearly, $u$ is well-defined and it holds 
    $$
    u \in L^{\infty}_{\loc}((0,T);L^2(\R^d)) \cap L^2_{\loc}((0,T);H^s_{\loc}(\R^d)) \cap L^1_{\loc}((0,T);L^1_{2s}(\R^d)).
    $$ 
    Moreover, since $p_t(\cdot,y) \to \delta_{\{y\}}$ in the sense of \eqref{eq:initial-mu} (see \autoref{lemma:heat-kernel-sol-property}), and therefore $u(t,\cdot) \to \mu$ in the sense of \eqref{eq:initial-mu}, as $t \searrow 0$. Hence, it only remains to verify that $u$ satisfies the solution property. To see this, we take $\psi \in C_c^{\infty}((0,T) \times \R^d)$. Then, we have the following computation, using the solution property for $(t,x) \mapsto p_t(x,y)$ for a.e. $y \in \R^d$ (see \autoref{lemma:heat-kernel-sol-property}), and Fubini's theorem:
    \begin{align}
    \label{eq:sol-check-convolution}
    \begin{split}
         \int_0^{\infty} \int_{\R^d} u(t,x) \partial_t \psi(t,x) \d x \d t &=  \int_0^{\infty} \int_{\R^d} \bigg(\int_{\R^d} p_t(x,y) \d \mu(y) \bigg) \partial_t \psi(t,x) \d x \d t \\
        &=  \int_{\R^d} \left( \int_0^{\infty}  \int_{\R^d} p_t(x,y) \partial_t \psi(t,x) \d x \d t \right) \d \mu(y) \\
        &=  \int_{\R^d} \bigg( \int_0^{\infty}   \cE^{(t)}(p_t(\cdot,y),\psi(t)) \d t \bigg) \d \mu(y) \\
        &=  \int_0^{\infty}   \cE^{(t)} \bigg( \int_{\R^d} p_t(\cdot,y) \d \mu(y),\psi(t) \bigg) \d t \\
        &=  \int_0^{\infty}   \cE^{(t)}( u(t),\psi(t)) \d t.
        \end{split}
    \end{align}
    In order to justify the  application of Fubini's theorem in the previous computation, we need to prove that the following two integrals converge: 
    \begin{align*}
       I_1 &:= \int_{\R^d} \int_0^{\infty} \int_{\R^d} |p_t(x,y)| |\partial_t \psi(t,x)| \d x \d t \d |\mu|(y), \\
       I_2 &:= \int_{\R^d} \int_0^{\infty}   |\cE^{(t)}(p_t(\cdot,y),\psi(t))| \d t \d |\mu|(y).
    \end{align*}
    For the integral $I_1$, we use that $\supp(\psi) \subset (T_1,T_2) \times B_R$ for some $0 < T_1 < T_2<T$ and $R>0$, and estimate by \eqref{eq:heat-kernel-upper}:
    \begin{align*}
        I_1 &\le c \int_{\R^d} \int_{T_1}^{T_2} \int_{B_R} |p_t(x,y)|  \d x \d t \d |\mu|(y) \\
        &\le c(T_1,T_2) \int_{\R^d} \left( \int_{B_R} \frac{\d x}{(1 + |x-y|)^{d+2s}}  \right) \d |\mu|(y) \\
        &\le c(T_1,T_2,R) \int_{\R^d}  \frac{\d |\mu|(y)}{(1 + |y|)^{d+2s}}  < \infty,
    \end{align*}
    where we used that $1+|y| \le c(1+|x-y|)$ for any $x \in B_R$ and $y \in \R^d$, for some constant $c > 0$, depending on $R$, and the assumption on $\mu$. 

    It remains to treat $I_2$. Again suppose $\supp(\psi) \subset (T_1,T_2) \times B_R$. We first estimate $\cE^{(t)}$ by splitting the domain of integration:
    \begin{align*}
        |\cE^{(t)}(p_t(\cdot,y),\psi(t))|&\le \int_{B_{2R}}\int_{B_{2R}}\frac{|p_t(x,y)-p_t(z,y)||\psi(t,x)-\psi(t,z)|}{|x-z|^{d+2s}}\d x\d z\\
        &\quad+2\int_{B_{R}}\int_{\R^d\setminus B_{2R}}\frac{|p_t(x,y)-p_t(z,y)||\psi(t,z)|}{|x-z|^{d+2s}}\d x\d z.
    \end{align*}
    Then, we plug this back into $I_2$ and apply H\"older's inequality, as follows:
    \begin{align*}
        I_2 &\le \Vert [\psi]_{H^s(B_{2R})} \Vert_{L^2(T_1,T_2)} \int_{\R^d} \Vert [p_t(\cdot,y)]_{H^s(B_{2R})} \Vert_{L^2(T_1,T_2)} \d |\mu|(y) \\
        &\quad + 2 \Vert \psi \Vert_{L^{\infty}((T_1,T_2) \times B_R)} \int_{\R^d} \int_{T_1}^{T_2}  \int_{B_R}\bigg(\int_{\R^d \setminus B_{2R}} \frac{|p_t(x,y) - p_t(z,y)|}{|x-z|^{d+2s} }\d z \bigg) \d x  \d t \d |\mu|(y) \\
        &=: I_{2,1} + I_{2,2}.
    \end{align*}
    Let us first prove $I_{2,1}$ is finite.
    Note that since $(t,x) \mapsto p_t(x,y)$ is a solution in $(T_1,T_2) \times \R^d$ for a.e. $y \in \R^d$ by \autoref{lemma:heat-kernel-sol-property}, we can apply the Caccioppoli inequality (see \cite[Proposition 2.1]{Lia24})  to deduce that for any $y \in \R^d$
    \begin{align*}
        \Vert [p_t(\cdot,y)]_{H^s(B_{2R})} \Vert^2_{L^2(T_1,T_2)} &\le c \int_{T_1}^{T_2}\int_{B_{4R}} p_t(x,y)^2 \d x \d t  \\
        &\qquad + c   \int_{T_1}^{T_2} \bigg(\int_{B_{4R}} p_t(x,y) \d x \bigg)\bigg( \int_{\R^d \setminus B_{4R}} \frac{p_t(x,y)}{|x|^{d+2s}} \d x \bigg) \d t,
        \end{align*}
    where $c=c(d,s,\lm,\Lm, R)$. The first term on the right-hand side is estimated by \eqref{eq:heat-kernel-upper} as
    \begin{align*}
        \int_{T_1}^{T_2}\int_{B_{4R}} p_t(x,y)^2 \d x \d t \le c \int_{B_{4R}} \frac{\d x}{(1 + |x-y|)^{2(d+2s)}}  \le \frac{c}{(1 + |y|)^{2(d+2s)}}
    \end{align*}
    for some $c=c(d,s,\lm,\Lm, R, T_1, T_2)$. Here, we also used that $(1 + |y|) \le c(R) (1 + |x-y|)$ for any $x \in B_{4R}$ and $y \in \R^d$ to obtain the last line. The second term is estimated by \eqref{eq:heat-kernel-upper} and \cite[Lemma~3.4]{LiWe24}:
        \begin{align*}
        & c  \bigg(\int_{T_1}^{T_2}\int_{B_{4R}} p_t(x,y) \d x \d t\bigg)   \cdot \bigg(\sup_{t \in (T_1,T_2)} \int_{\R^d \setminus B_{4R}} \frac{p_t(x,y)}{|x|^{d+2s}} \d x \d t\bigg) \\
        &\qquad\le c \bigg(\int_{B_{4R}} \frac{\d x}{(1 + |x-y|)^{d+2s}} \bigg)  \cdot \bigg( \int_{\R^d \setminus B_{4R}} \frac1{(1 + |x-y|)^{d+2s}}\frac{\d x}{|x|^{d+2s}} \bigg)  \\
        &\qquad\le \frac{c}{(1 + |y|)^{2(d+2s)}},
    \end{align*}
    for some $c=c(d,s,\lm,\Lm, R, T_1, T_2)$. Hence, we have shown
    \begin{align*}
        I_{2,1} \le c\,\Vert [\psi]_{H^s(B_{2R})} \Vert_{L^2(T_1,T_2)} \int_{\R^d}\frac{\d |\mu|(y)}{(1 + |y|)^{d+2s}}  < \infty.
    \end{align*}
For $I_{2,2}$, we observe, using again \eqref{eq:heat-kernel-upper}
\begin{align*}
    \int_{T_1}^{T_2} \int_{B_R} |p_t(x,y)| \bigg( \int_{\R^d \setminus B_{2R}} \frac{\d z}{|x-z|^{d+2s}}  \bigg) \d x \d t 
    &\le c \int_{B_R} \frac{\d x}{(1 + |x-y|)^{d+2s}}  \\
    &\le  \frac{c}{(1 + |y|)^{d+2s}}, 
\end{align*}
for some $c=c(d,s,\lm,\Lm, R, T_1, T_2)$. Moreover, using that $(1+|z-y|) \le c (1 + |x+z-y|)$ for any $x \in B_R$, $y \in \R^d$ and $z \in \R^d \setminus B_{2R}$:
\begin{align*}
    \int_{T_1}^{T_2} \int_{B_R} &\bigg( \int_{\R^d \setminus B_{2R}} \frac{|p_t(z,y)|}{|x-z|^{d+2s}} \d z \bigg) \d x \d t \\
    &\le c \sup_{x \in B_R} \int_{\R^d \setminus B_{R}(x)} \frac1{(1 + |z-y|)^{d+2s}} \frac{\d z}{|x-z|^{d+2s}}  \\
    &\le c \sup_{x \in B_R} \int_{\R^d \setminus B_{R}} \frac1{(1 + |x+z-y|)^{d+2s}} \frac{\d z}{|z|^{d+2s}}  \\
    &\le c \int_{\R^d \setminus B_{R}} \frac1{(1 + |z-y|)^{d+2s}} \frac{\d z}{|z|^{d+2s}} \\
    &\le \frac{c}{(1 + |y|)^{d+2s}},
\end{align*}
    for some $c=c(d,s,\lm,\Lm, R, T_1, T_2)$, where we used \cite[Lemma~3.4]{LiWe24} in the last step.
    Thus, we have shown
    \begin{align*}
        I_{2,2} \le  c\, \Vert \psi \Vert_{L^{\infty}((T_1,T_2) \times B_R)} \int_{\R^d} \frac{\d |\mu|(y)}{(1 + |y|)^{d+2s}}  < \infty,
    \end{align*}
which implies $I_2 < \infty$.

    Altogether, the applications of Fubini's theorem in \eqref{eq:sol-check-convolution} are justified, and therefore, we have shown that for any $\psi \in C_c^{\infty}((0,T) \times \R^d)$ it holds
    \begin{align*}
        - \int_0^{\infty} \int_{\R^d} u(t,x) \partial_t \psi(t,x) \d x \d t  +\int_0^{\infty}  \cE^{(t)}( u(t),\psi(t)) \d t = 0. 
    \end{align*}
    By a standard density argument, we conclude that $u$ has the solution property described in \autoref{def:global-sol}, as desired. The proof is complete.
\end{proof}

\subsection{Proof of the Widder-type theorem}
We are now in a position to give the proof of the Widder-type theorem

\begin{proof}[Proof of \autoref{thm:widder}]
    First, because of $u \in L^1((0,T) ; L^1(\R^d;w))$, there exists $T_o \in (0,T)$ such that $u(T_o,\cdot) \in L^1(\R^d;w)$. Thus, by \autoref{lemma:ex-un-initial-trace}  there exists a unique $\mu \in \cM_+$ satisfying the growth condition
\begin{align*}
        \int_{\R^d} \frac{\d \mu(x)}{(1 + |x|)^{d+2s}}  \le c\, \Vert u(T_o) \Vert_{L^1(\R^d;w)},
    \end{align*}
such that $u(0) = \mu$. Therefore, $u$ is a solution to the Cauchy problem with initial datum $\mu$ in the sense of \autoref{def:Cauchy}. On the other hand, by \autoref{lemma:existence}, we have that
\begin{align*}
    v(t,x) := \int_{\R^d} p_t(x,y) \d \mu(y)
\end{align*}
is also a solution to the Cauchy problem with initial datum $\mu$. Finally, by \autoref{lemma:uniqueness}, we deduce that $u = v$, and the proof is complete. 
\end{proof}

\section{Harnack-type estimates for global solutions}
\label{sec:est-global}

We are now in a position to prove \autoref{thm:est-global} and \autoref{thm:est-global-}.

\begin{proof}[Proof of \autoref{thm:est-global} and \autoref{thm:est-global-}]
    Up to resorting to the new function 
    \begin{align}
        \label{eq:v}
        v(t,x):=u(\tau t, \tau^{1/2s}(x+z))
    \end{align}
    with $\tau>0$ and $z \in \R^d$, it suffices to consider a global solution $u$ defined in $(0,1]\times\R^d$ and to show the desired estimates for $t=1$ and $y=0$. To this end, we rely on the elliptic-type Harnack estimate established in \cite[Theorem~1.2]{LiWe24}. 
     For $\xi\in(0,1]$ and $x_o = x$ in the supremum estimate of \cite[Theorem~1.2]{LiWe24} we obtain that
    \begin{align*}
        u(\xi,x)&\le c \,\xi \int_{\R^d}\frac{u(\xi,y)}{(\xi^{\frac{1}{2s}}+|x-y|)^{d+2s}}\d y\\
        &\le c\, \xi \sup_{z \in \R^d} \bigg( \frac{1 + |z|}{\xi^{\frac{1}{2s}} + |x-z|} \bigg)^{d+2s} \int_{\R^d} \frac{u(\xi,y)}{(1 + |y|)^{d+2s}}\d y
    \end{align*}
    for some $c=c(d,s,\lambda,\Lambda)$. To proceed, observe that for any $x,\,z\in\R^d$ we have the simple algebraic fact, using that $\xi \le 1$:
    \begin{align}
    \label{eq:algebra}
    \frac{1 + |z|}{\xi^{\frac{1}{2s}} + |x-z|} \le \frac{1 + |x-z|}{\xi^{\frac{1}{2s}} + |x-z|} + \frac{|x|}{\xi^{\frac{1}{2s}} + |x-z|} \le \xi^{-{\frac{1}{2s}}} +  \xi^{-\frac{1}{2s}}|x|.
    \end{align}
    Moreover, we first apply the time-insensitive $L^1(\R^d;w)$ estimate from \cite[Proposition~6.1]{LiWe24} and then apply
     the infimum estimate of \cite[Theorem~1.2]{LiWe24} taking $x_o=0$ and $\tau=1$. In this way, we obtain that
    \begin{align*}
        \int_{\R^d} \frac{u(\xi,y)}{(1 + |y|)^{d+2s}} \d y &\le c \int_{\R^d} \frac{u(1,y)}{(1 + |y|)^{d+2s}} \d y \le c\, u(1,0)
    \end{align*}
    for some $c=c(d,s,\lambda,\Lambda)$.
    Combining these observations yields 
    \begin{align*}
        u(\xi,x) &\le c \, \xi \big(\xi^{-{\frac{1}{2s}}} + \xi^{-\frac{1}{2s}}|x| \big)^{d+2s} u(1,0) \\
        &\le c \, \xi^{-\frac{d}{2s}} (1+|x|)^{d+2s}  u(1,0), 
    \end{align*}
    which immediately implies the first desired conclusion.
    As for the second conclusion, we use the representation formula in \autoref{thm:widder} and the upper estimate in \autoref{thm:hkb} to obtain that
    \begin{align*}
        u(1,x)&\le c\int_{\R^d}\frac{\d\mu(y)}{(1+|x-y|)^{d+2s}} \\
        &\le c\, (1+|x|)^{d+2s} \int_{\R^d} \frac{\d\mu(y)}{(1 + |y|)^{d+2s}}.
    \end{align*}
    Here, again we used the above algebraic fact \eqref{eq:algebra} and apply \eqref{eq:v} with $z = 0$, which concludes the proof.
\end{proof}

\appendix

\section{Uniqueness for signed solutions}
\label{sec:signed-sol}

We include a proof of uniqueness of signed solutions which satisfy a certain additional regularity at time zero, namely
 \begin{align}
 \label{eq:u-nice-up-to-zero}
     u \in L^{\infty}((0,T) ; L^2_{\loc}(\R^d)) \cap L^2((0,T) ; H^s_{\loc}(\R^d)) \cap L^1((0,T); L^1(\R^d;w)).
\end{align}

The following result is in analogy to \cite[Theorem 2, p.639]{Aro68}. 

\begin{proposition}
\label{lemma:uniqueness-signed}
        Assume that the kernel $K$ satisfies the upper bound \eqref{eq:up-bd} and the coercivity condition \eqref{eq:coercive}. Let $u_1$ and $u_2$ be, possibly signed, weak solutions to the Cauchy problem in $(0,T] \times \R^d$ in the sense of \autoref{def:Cauchy} with $u_1(0) = \mu = u_2(0)$ for some $\mu \in \cM_+$. In addition, assume that $u_1$ and $u_2$ satisfy \eqref{eq:u-nice-up-to-zero}. Then, we have $u_1 = u_2$ in $(0,T] \times \R^d$.
\end{proposition}

As the first step we show that solutions satisfying \eqref{eq:u-nice-up-to-zero} can be extended to negative times if the initial datum is zero.

\begin{lemma}
\label{lemma:ext-neg-times}
Let $u$ be a local weak solution in $(0,T) \times \Omega$ in the sense of \autoref{def:local-sol} with $u(0) = 0$ in the sense of \eqref{eq:initial-mu}, assuming in addition that
\begin{align}
\label{eq:u-nice-up-to-zero-loc}
    u \in L^{\infty}((0,T) ; L^2_{\loc}(\Omega)) \cap L^2((0,T) ; H^s_{\loc}(\Omega)) \cap L^1((0,T); L^1(\R^d;w)).
\end{align}
Define $\bar{u}$ such that $\bar{u}(t) \equiv u(t)$ for $t > 0$ and $\bar{u}(t) \equiv 0$ for $t \le 0$; define $\bar{K}$ such that $\bar{K}(t) \equiv K(t)$ for $t > 0$ and $\bar{K}(t) \equiv |x-y|^{-d-2s}$ for $t \le 0$. Then, $\bar{u}$ is a local weak solution in $(-\infty,T) \times \Omega$ in the sense of \autoref{def:local-sol} with respect to $\bar{\mathcal{L}}_t$ driven by the kernel $\bar{K}$.
\end{lemma}

\begin{proof}
First, we claim that for any $\phi \in C^{\infty}_c((-\infty,T) \times \Omega)$ it holds
\begin{align}
\label{eq:weakL2-initial}
    \lim_{t \searrow 0} \int_{\Omega} u(t,x) \phi(t,x) \d x = 0.
\end{align}
In fact, for any $t > 0$ we consider
\begin{align*}
    \int_{\Omega} u(t,x) \phi(t,x) \d x = \int_{\Omega} u(t,x) [\phi(t,x) - \phi(0,x) ]\d x + \int_{\Omega} u(t,x) \phi(0,x) \d x.
\end{align*}
 As $t \searrow 0$, the first term on the left-hand side converges to zero since $u \in L^{\infty}((0,T);L^2_{\loc}(\Omega))$ and
\begin{align*}
    \left|\int_{\Omega} u(t,x) [\phi(t,x) - \phi(0,x) ]\d x \right| \le \sup_{t\in(0,T)}\Vert u \Vert_{L^2(\supp(\phi(t,\cdot)))} \|\phi(t,\cdot) - \phi(0,\cdot)\|_{L^2(\Om)},
\end{align*}
whereas the second term also converges to zero thanks to \eqref{eq:initial-mu}.
Therefore,  we have \eqref{eq:weakL2-initial}, as claimed.

Moreover, consider any $0 < t_1 < t_2 < T$ and introduce $\zeta_k \in C_c^1(0,T)$ such that $\zeta_k = 1$ in $[t_1,t_2]$, $\zeta_k = 0$ in $(0,t_1 - \frac{1}{k}] \cup [t_2+ \frac{1}{k},T)$, and linearly interpolated otherwise. Test the weak formulation for $u$ with $\zeta_k \phi$, and then send $k \to \infty$ to get
\begin{align}
\label{eq:t1t2}
    \int_{\Omega} u(t,x) \phi(t,x) \d x\bigg|_{t_1}^{t_2} + \int_{t_1}^{t_2} \int_{\Omega}  -\partial_t \phi(t,x) u(t,x) \d x \d t + \int_{t_1}^{t_2}\cE^{(t)}(u(t),\phi(t)) \d t = 0.
\end{align}

Therefore, sending $t_1 \searrow 0$ in \eqref{eq:weakL2-initial} with the aid of \eqref{eq:t1t2}, we deduce
\begin{align*}
    \int_{\Omega} u(t_2,x) \phi(t_2,x) \d x + \int_{0}^{t_2} \int_{\Omega} -\partial_t \phi(t,x) u(t,x) \d x \d t + \int_{0}^{t_2}\cE^{(t)}(u(t),\phi(t)) \d t = 0.
\end{align*}
Note that our assumptions on $u$ and $\phi$ allow us to take limits by dominated convergence, where we are using that $u$ is well-behaved up to $t=0$, see \eqref{eq:u-nice-up-to-zero-loc}. 
In particular, upon sending $t_2 \nearrow T$ and using that $\phi(T) \equiv 0$ by assumption, we deduce
\begin{align}
\label{eq:positive-times}
    \int_{0}^{T} \int_{\Omega} -\partial_t \phi(t,x) u(t,x) \d x \d t+ \int_{0}^{T}\cE^{(t)}(u(t),\phi(t)) \d t = 0.
\end{align}

On the other hand, observe that by definition of extension we have
\begin{align}
\label{eq:negative-times}
    \int^{0}_{-\infty} \int_{\Omega} -\partial_t \phi(t,x) \bar{u}(t,x) \d x \d t + \int^{0}_{-\infty}\bar{\cE}^{(t)}(\bar{u}(t),\phi(t)) \d t = 0.
\end{align}
Altogether, summing up \eqref{eq:positive-times} and \eqref{eq:negative-times}, this yields 
\begin{align*}
    \int^{T}_{-\infty} \int_{\Omega} -\partial_t \phi(t,x) \bar{u}(t,x) \d x \d t + \int^{T}_{-\infty}\bar{\cE}^{(t)}(\bar{u}(t),\phi(t)) \d t = 0,
\end{align*}
for any $\phi \in C^{\infty}_c((-\infty,T) \times \Omega)$. By density, this implies that $\bar u$ is a local weak solution in the sense of \autoref{def:local-sol} with respect to $\bar{\mathcal{L}}_t$ driven by the kernel $\bar{K}$.
\end{proof}

With \autoref{lemma:ext-neg-times} at hand, we can prove \autoref{lemma:uniqueness-signed}. Essentially, this follows from the Liouville type property which can be deduced from H\"older regularity of weak solutions (see for instance \cite[Remark~1.1]{Lia24b} and \cite[Theorem 2.1]{FeRo17}).

\begin{proof}[Proof of \autoref{lemma:uniqueness-signed}]
    Let us take $u = u_1 - u_2$ and observe that $u$ solves the Cauchy problem in $(0,T) \times \R^d$ with $u(0) = 0$ and $u$ satisfies \eqref{eq:u-nice-up-to-zero}. In particular, by \autoref{lemma:ext-neg-times}, we can extend $u$ to a solution $\bar{u}$ in $(-\infty,T) \times \R^d$ by extending the kernel $K$ also to negative times.
    
    Next, we show that $u$ is bounded in $(-\infty,T) \times \R^d$. To see this, we apply the local boundedness estimate
\[
\sup_{(T-R^{2s},T]\times B_{R}} |\bar u| \le  c \int_{T-2 R^{2s}}^{T} \int_{\R^d} \frac{|\bar u(t,x)|}{(R+|x|)^{d+2s} }\d x \d t,
\]
where $R>0$ is arbitrary and $c=c(d,s,\lm,\Lm)$. This result follows from \cite[Proposition 3.6]{LiWe24}, and also holds true for signed solutions, first taking $t_o=1$ and $x_o=0$, and then applying a proper scaling argument. By sending $R \nearrow \infty$ and using that $\bar{u}(t,\cdot) = 0$ for $t \le 0$, and $u \in L^1((0,T);L^1(\R^d;w))$, we obtain
\begin{align*}
    \sup_{(-\infty,T]\times \R^d} |\bar u| \le  c \int_{0}^{T} \int_{\R^d} \frac{|u(t,x)|}{(1+|x|)^{d+2s} }\d x \d t < \infty,
\end{align*}
so $u$ is indeed globally bounded. 

From here, we can deduce that $u \equiv 0$ by the Liouville theorem (see also \cite[Remark 1.1]{Lia24b}). Indeed, using the H\"older regularity estimate from \cite[Theorem 3.8]{LiWe24}, we deduce that there exist $c>1$ and $\alpha \in (0,1)$ depending on $\{d,s,\lm,\Lm\}$, such that for any $R > 0$ we have
\begin{align*}
    [\bar{u}]_{C^{\alpha}((T-R^{2s},T] \times B_R)} \le c\, R^{-\alpha} \sup_{(-\infty,T]\times \R^d} |\bar u|.
\end{align*}
Hence, sending $R \nearrow \infty$, we deduce that
\begin{align*}
    [\bar{u}]_{C^{\alpha}((-\infty,T) \times \R^d)} = 0,
\end{align*}
and therefore, $\bar{u}$ must be constant. Since $\bar{u}(0) = 0$, we deduce that $u \equiv 0$, as desired.
\end{proof}

\section{Continuity of solutions up to the initial level}
\label{sec:cont-initial}

The purpose of this section is to show the local continuity of weak solutions at the initial level. This is an essential ingredient in the proof of \autoref{lemma:heat-kernel-sol-property} as it allows us to verify that $p_{\eta,t}(\cdot,y) \to \delta_{\{y\}}$ in the sense of \eqref{eq:initial-mu} for any $y \in \R^d$.

\begin{proposition}\label{prop:initial}
    Assume that the kernel $K$ satisfies the upper bound \eqref{eq:up-bd} and the coercivity condition \eqref{eq:coercive}. Let $u$ be a local weak solution to \eqref{eq:PDE-domains} in $(0,T]\times \Om$ in the sense of \autoref{def:local-sol}. Suppose that $|u|\le \frac{1}{2}M$ in $(0,T]\times\R^d$ for some $M>0$ and that $u(t,\cdot)\to u_o\in C_{\loc}(\Om)$ as $t\searrow0$ in the sense of $L^2_{\loc}(\Om)$. Then $u$ is locally continuous in $[0,T]\times \Om$.
\end{proposition}

\begin{proof}
It suffices to show that given a point $x_o\in\Om$, $u$ is continuous at $(0,x_o)$.
    In order to prove it, consider the ``forward" cylinder $Q_R:=[0,R^{2s}]\times B_{R}(x_o)$ inside $[0,T]\times \Om$ and the quantities
    \[
    \mu^+:=\sup_{Q_R}u,\qquad \mu^-:=\inf_{Q_R} u.
    \]
    For ease of notation, we omit $x_o$ from $B_R(x_o)$. Note that by $|u| \le \frac{1}{2}M$ we have $\osc_{Q_R} u \le M$. We will now establish oscillation decay in smaller cylinders. To do so, we suppose that one of the following holds true:
    \begin{equation}\label{alter:0}
    \mu^+-\sup_{B_R}u_o\ge\tfrac14 M\qquad\text{or}\qquad\inf_{B_R}u_o -\mu^-\ge\tfrac14 M.
    \end{equation}
    Otherwise we end up with
    \begin{equation}\label{om-u0}
    \osc_{Q_R} u\le \tfrac12 M + \osc_{B_R}u_o.
    \end{equation}
    Let us assume that the first property in \eqref{alter:0} holds true. According to \cite[Lemma~3.2]{Lia24b} (which remains true for kernels satisfying the coercivity condition \eqref{eq:coercive} instead of \eqref{eq:low-bd}) with $\xi\om=\frac14M$, $\varrho=\dl R$, and $\mathcal{Q}=Q_R$, we find $\nu\in(0,1)$ depending on $\{d,s,\lm,\Lm\}$ such that either 
    \[
    \dl^{2s}{\rm Tail}[(u-\mu^+)_+; Q_R]>\tfrac14M,
    \]
    where we define
    \begin{align}
        \tail(v;Q_R) = R^{2s} \sup_{t \in (0,R^{2s})} \int_{\R^d \setminus B_R(x_o)} \frac{|v(t,x)|}{|x-x_o|^{d+2s}} \d x,
    \end{align}
    or it holds
    \begin{equation}\label{osc:0}
    \mu^+- u\ge\tfrac18M \quad \text{a.e. in}\> [0,\nu(\dl R)^{2s}]\times B_{\dl R}(x_o).
    \end{equation}
    Setting $\sig:=\dl\nu^{\frac1{2s}}$, \eqref{osc:0} gives
    \begin{equation*}
        \osc_{Q_{\sig R}}u\le \mu^+ - \tfrac{1}{8}M - \inf_{Q_{\sig R}} u \le \tfrac78 M,
    \end{equation*}
    whereas the tail alternative will not occur because $${\rm Tail}[(u-\mu^+)_+; Q_R]\le c\,M$$ for some $c=c(d,s,\lm,\Lm)$ and we can choose $\dl$ to be so small that $\dl^{2s}c\le \frac{1}{4}$. If the second case of \eqref{alter:0} holds true, we will end up with the same oscillation estimate as above. Therefore, denoting $R_1:= \sig R$ and combining \eqref{om-u0} and \eqref{osc:0} we arrive at
    \begin{equation*}
        \osc_{Q_{R_1}} u\le \max\Big\{\tfrac78 M,\, \tfrac12 M + \osc_{B_R}u_o\Big\}=:M_1.
    \end{equation*}

    Now, suppose we have constructed the following up to $i=1,\dots, j$:
    \begin{equation*}
    \left\{
        \begin{array}{cc}
            \displaystyle R_o=R,\quad R_i=\sig R_{i-1},\\
            \displaystyle M_o=M,\quad M_i= \max\Big\{\tfrac78 M_{i-1},\, \tfrac12 M_{i-1} + \osc_{B_{R_{i-1}}}u_o\Big\}, \\
             \displaystyle \mu_i^+=\sup_{Q_{R_i}}u, \quad \mu_i^-=\inf_{Q_{R_i}}u,\,\quad \osc_{Q_{R_i}} u\le M_i .
        \end{array}\right.
    \end{equation*}
    To implement induction we show that it continues to hold for $i=j+1$. Indeed, we consider the alternative
        \begin{equation*}
    \mu_j^+-\sup_{B_{R_j}}u_o\ge\tfrac14 M_j\qquad\text{or}\qquad\inf_{B_{R_j}}u_o -\mu_j^-\ge\tfrac14 M_j;
    \end{equation*}
    otherwise we end up with
    \begin{equation*}
    \osc_{Q_{R_j}} u\le \tfrac12 M_j + \osc_{B_{R_j}}u_o.
    \end{equation*}
    Repeating the same reasoning as before will yield 
    \begin{equation*}
        \osc_{Q_{\sig R_j}} u\le \max\Big\{\tfrac78 M_j,\, \tfrac12 M_j + \osc_{B_{R_j}}u_o\Big\}=:M_{j+1}.
    \end{equation*}
    A key for this oscillation estimate to hold is the tail estimate $${\rm Tail}[(u-\mu_j^+)_+; Q_{R_j}]\le c\,M_j$$ for some $c=c(d,s,\lm,\Lm)$. This estimate can be proven as \cite[(4.9)]{Lia24b} noticing that $M_j\ge(\frac78)^{j-i+1}M_{i-1}$ for any $j\ge i\ge1$.
    Thus, the induction is completed. 
    
    To obtain the continuity of $u$ at $(0,x_o)$ we argue as follows. Given $\varepsilon>0$, we first fix $R$ satisfying
    \[
    \osc_{B_{R}}u_o\le \varep.
    \]
    Then, by construction of $M_i$ we estimate
    \[
    M_{i}\le \tfrac78 M_{i-1} +\varep,\quad\forall\, i\in \N
    \]
    and iterate to get
    \[
    M_{i}\le (\tfrac78)^i M+ \eps \sum_{j=0}^i \left( \tfrac{7}{8} \right)^j \le (\tfrac78)^i M+8\varep,\quad\forall\, i\in \N.
    \]
    Finally, we choose $i$ so large that $M_i<9\varepsilon$ and combine it with the oscillation estimate to conclude that
    \[
    \osc_{Q_{R_i}}u\le M_i \le 9\varepsilon.
    \]
    The continuity at $(0,x_o)$ is thus proven.
\end{proof}

\begin{remark}
    It is possible to extract an explicit modulus of continuity out of the above argument. Moreover, the global boundedness of $u$ can be relaxed as in \cite{Lia24b}.
\end{remark}


\begin{thebibliography}{dPQRV12}

\bibitem[APT22]{APT22}
K.~Adimurthi, H.~Prasad, and V.~Tewary.
\newblock Local {H}{\"o}lder regularity for nonlocal parabolic $p$-{L}aplace
  equations.
\newblock {\em arXiv:2203.13082 (to appear in Ann. Sc. Norm. Super. Pisa Cl.
  Sci. (5))}, 2022.

\bibitem[Aro68]{Aro68}
D.~G. Aronson.
\newblock Non-negative solutions of linear parabolic equations.
\newblock {\em Ann. Scuola Norm. Sup. Pisa Cl. Sci. (3)}, 22:607--694, 1968.

\bibitem[Aro71]{Aro71}
D.~G. Aronson.
\newblock Addendum: ``{N}on-negative solutions of linear parabolic equations''\
  ({A}nn. {S}cuola {N}orm. {S}up. {P}isa (3) {\bf 22} (1968), 607--694).
\newblock {\em Ann. Scuola Norm. Sup. Pisa Cl. Sci. (3)}, 25:221--228, 1971.

\bibitem[BGK09]{BGK09}
M.~Barlow, A.~Grigor'yan, and T.~Kumagai.
\newblock Heat kernel upper bounds for jump processes and the first exit time.
\newblock {\em J. Reine Angew. Math.}, 626:135--157, 2009.

\bibitem[BPSV14]{BPSV14}
B.~Barrios, I.~Peral, F.~Soria, and E.~Valdinoci.
\newblock A {W}idder's type theorem for the heat equation with nonlocal
  diffusion.
\newblock {\em Arch. Ration. Mech. Anal.}, 213(2):629--650, 2014.

\bibitem[BL02]{BaLe02}
R.~Bass and D.~Levin.
\newblock Transition probabilities for symmetric jump processes.
\newblock {\em Trans. Amer. Math. Soc.}, 354(7):2933--2953, 2002.

\bibitem[BG60]{BlGe60}
R.~Blumenthal and R.~Getoor.
\newblock Some theorems on stable processes.
\newblock {\em Trans. Amer. Math. Soc.}, 95:263--273, 1960.

\bibitem[BSV17]{BSV17}
M.~Bonforte, Y.~Sire, and J.~V\'{a}zquez.
\newblock Optimal existence and uniqueness theory for the fractional heat
  equation.
\newblock {\em Nonlinear Anal.}, 153:142--168, 2017.


\bibitem[BKO23]{BKO23}
S.-S. Byun, H.~Kim, and J.~Ok.
\newblock Local {H}\"older continuity for fractional nonlocal equations with
  general growth.
\newblock {\em Math. Ann.}, 387(1-2):807--846, 2023.

\bibitem[CCV11]{CCV11}
L.~Caffarelli, C.~H. Chan, and A.~Vasseur.
\newblock Regularity theory for parabolic nonlinear integral operators.
\newblock {\em J. Amer. Math. Soc.}, 24(3):849--869, 2011.


\bibitem[CKS87]{CKS87}
E.~A. Carlen, S.~Kusuoka, and D.~W. Stroock.
\newblock Upper bounds for symmetric {M}arkov transition functions.
\newblock {\em Ann. Inst. H. Poincar\'e{} Probab. Statist.}, 23(2):245--287, 1987.
  
\bibitem[CKW22a]{CKW22b}
J.~Chaker, M.~Kim, and M.~Weidner.
\newblock Regularity for nonlocal problems with non-standard growth.
\newblock {\em Calc. Var. Partial Differential Equations}, 61(6):Paper No. 227, 31, 2022.
  
\bibitem[CKW23]{CKW23}
J.~Chaker, M.~Kim, and M.~Weidner.
\newblock Harnack inequality for nonlocal problems with non-standard growth.
\newblock {\em Math. Ann.}, 386(1-2):533--550, 2023.
  
\bibitem[CKKW21]{CKKW21}
Z.-Q. Chen, P.~Kim, T.~Kumagai, and J.~Wang.
\newblock Heat kernel upper bounds for symmetric {M}arkov semigroups.
\newblock {\em J. Funct. Anal.}, 281(4):Paper No. 109074, 40, 2021.

\bibitem[CK03]{ChKu03}
Z.-Q. Chen and T.~Kumagai.
\newblock Heat kernel estimates for stable-like processes on {$d$}-sets.
\newblock {\em Stochastic Process. Appl.}, 108(1):27--62, 2003.

\bibitem[CK08]{ChKu08}
Z.-Q. Chen and T.~Kumagai.
\newblock Heat kernel estimates for jump processes of mixed types on metric
  measure spaces.
\newblock {\em Probab. Theory Related Fields}, 140(1-2):277--317, 2008.

\bibitem[CKW21]{CKW21}
Z.-Q. Chen, T.~Kumagai, and J.~Wang.
\newblock Stability of heat kernel estimates for symmetric non-local
  {D}irichlet forms.
\newblock {\em Mem. Amer. Math. Soc.}, 271(1330):v+89, 2021.

\bibitem[CKW22b]{CKW22c}
Z.-Q. Chen, T.~Kumagai, and J.~Wang.
\newblock Heat kernel estimates for general symmetric pure jump {D}irichlet
  forms.
\newblock {\em Ann. Sc. Norm. Super. Pisa Cl. Sci. (5)}, 23(3):1091--1140, 2022.

\bibitem[Coz17]{Coz17}
M.~Cozzi.
\newblock Regularity results and {H}arnack inequalities for minimizers and
  solutions of nonlocal problems: a unified approach via fractional {D}e
  {G}iorgi classes.
\newblock {\em J. Funct. Anal.}, 272(11):4762--4837, 2017.

\bibitem[dTEJ17]{DEJ17}
F.~del Teso, J.~Endal, and E.~Jakobsen.
\newblock Uniqueness and properties of distributional solutions of nonlocal
  equations of porous medium type.
\newblock {\em Adv. Math.}, 305:78--143, 2017.

\bibitem[DS23]{DeSi23}
M.~Dembny and M.~Sierzkega.
\newblock A sharp {H}arnack bound for a nonlocal heat equation.
\newblock {\em arXiv:2303.08186}, 2023.

\bibitem[dPQRV12]{DQRV12}
A.~de~Pablo, F.~Quir\'os, A.~Rodr\'iguez, and J.~L. V\'azquez.
\newblock A general fractional porous medium equation.
\newblock {\em Comm. Pure Appl. Math.}, 65(9):1242--1284, 2012.

\bibitem[DCKP14]{DKP14}
A.~Di~Castro, T.~Kuusi, and G.~Palatucci.
\newblock Nonlocal {H}arnack inequalities.
\newblock {\em J. Funct. Anal.}, 267(6):1807--1836, 2014.

\bibitem[DCKP16]{DKP16}
A.~Di~Castro, T.~Kuusi, and G.~Palatucci.
\newblock Local behavior of fractional {$p$}-minimizers.
\newblock {\em Ann. Inst. H. Poincar\'{e} C Anal. Non Lin\'{e}aire},
  33(5):1279--1299, 2016.
  
\bibitem[EG15]{EvGa15}
L.~Evans and R.~Gariepy.
\newblock {\em Measure theory and fine properties of functions}.
\newblock Textbooks in Mathematics. CRC Press, Boca Raton, FL, revised edition,
  2015.

\bibitem[FK13]{FeKa13}
M.~Felsinger and M.~Kassmann.
\newblock Local regularity for parabolic nonlocal operators.
\newblock {\em Comm. Partial Differential Equations}, 38(9):1539--1573, 2013.

\bibitem[FRRO17]{FeRo17}
X.~Fern\'andez-Real and X.~Ros-Oton.
\newblock Regularity theory for general stable operators: parabolic equations.
\newblock {\em J. Funct. Anal.}, 272(10):4165--4221, 2017.

\bibitem[GHH17]{GHH17}
A.~Grigor'yan, E.~Hu, and J.~Hu.
\newblock Lower estimates of heat kernels for non-local {D}irichlet forms on
  metric measure spaces.
\newblock {\em J. Funct. Anal.}, 272(8):3311--3346, 2017.

\bibitem[GHH18]{GHH18}
A.~Grigor'yan, E.~Hu, and J.~Hu.
\newblock Two-sided estimates of heat kernels of jump type {D}irichlet forms.
\newblock {\em Adv. Math.}, 330:433--515, 2018.

\bibitem[GHH23]{GHH23}
A.~Grigor'yan, E.~Hu, and J.~Hu.
\newblock Off-diagonal lower estimates and {H}\"older regularity of the heat
  kernel.
\newblock {\em Asian J. Math.}, 27(5):675--770, 2023.

\bibitem[GHH24]{GHH24}
A.~Grigor'yan, E.~Hu, and J.~Hu.
\newblock Tail estimates and off-diagonal upper bounds of the heat kernel.
\newblock {\em preprint, https://www.math.uni-bielefeld.de/~grigor/tp-ueq.pdf},
  2024.

\bibitem[GHL14]{GHL14}
A.~Grigor'yan, J.~Hu, and K.-S. Lau.
\newblock Estimates of heat kernels for non-local regular {D}irichlet forms.
\newblock {\em Trans. Amer. Math. Soc.}, 366(12):6397--6441, 2014.

\bibitem[GMP15]{GMP15}
G.~Grillo, M.~Muratori, and F.~Punzo.
\newblock Fractional porous media equations: existence and uniqueness of weak
  solutions with measure data.
\newblock {\em Calc. Var. Partial Differential Equations}, 54(3):3303--3335,
  2015.

\bibitem[GQSV25]{GQSV25}
I.~Gonz{\'a}lvez, F.~Quir{\'o}s, F.~Soria, and Z.~Vondracek.
\newblock On the nonlocal heat equation for certain {L}{\'e}vy operators and
  the uniqueness of positive solutions.
\newblock {\em arXiv:2504.04246}, 2025.

\bibitem[Kan23]{Kan23}
J.~Kang.
\newblock Heat kernel estimates for symmetric jump processes with anisotropic
  jumping kernels.
\newblock {\em Proc. Amer. Math. Soc.}, 151(1):385--399, 2023.

\bibitem[Kas09]{Kas09}
M.~Kassmann.
\newblock A priori estimates for integro-differential operators with measurable
  kernels.
\newblock {\em Calc. Var. Partial Differential Equations}, 34(1):1--21, 2009.

\bibitem[KW23]{KaWe23}
M.~Kassmann and M.~Weidner.
\newblock Upper heat kernel estimates for nonlocal operators via {A}ronson's
  method.
\newblock {\em Calc. Var. Partial Differential Equations}, 62(2):Paper No. 68,
  27, 2023.

\bibitem[KW24]{KaWe23b}
M.~Kassmann and M.~Weidner.
\newblock The parabolic {H}arnack inequality for nonlocal equations.
\newblock {\em Duke Math. J.}, 173(17):3413--3451, 2024.

\bibitem[Krz64]{Krz64}
M.~Krzy\.za\'nski.
\newblock Sur les solutions non n\'egatives de l'\'equation lin\'eaire normale
  parabolique.
\newblock {\em Rev. Roumaine Math. Pures Appl.}, 9:393--408, 1964.

\bibitem[Lia24a]{Lia24}
N.~Liao.
\newblock {H}arnack estimates for nonlocal drift-diffusion equations.
\newblock {\em arXiv:2402.11986v2}, 2024.

\bibitem[Lia24b]{Lia24b}
N.~Liao.
\newblock H\"{o}lder regularity for parabolic fractional {$p$}-{L}aplacian.
\newblock {\em Calc. Var. Partial Differential Equations}, 63(1):Paper No. 22,
  34 pp., 2024.

\bibitem[Lia24c]{Lia24c}
N.~Liao.
\newblock On the modulus of continuity of solutions to nonlocal parabolic
  equations.
\newblock {\em J. Lond. Math. Soc.}, 110(3):Paper No. e12985, 30 pp., 2024.

\bibitem[LW24]{LiWe24}
N.~Liao and M.~Weidner.
\newblock Time-insensitive nonlocal parabolic {H}arnack estimates.
\newblock {\em arXiv:2409.20097 (to appear in Proc. Lond. Math. Soc.)}, 2024.

\bibitem[MM13a]{MaMi13a}
Y.~Maekawa and H.~Miura.
\newblock On fundamental solutions for non-local parabolic equations with
  divergence free drift.
\newblock {\em Adv. Math.}, 247:123--191, 2013.

\bibitem[MM13b]{MaMi13b}
Y.~Maekawa and H.~Miura.
\newblock Upper bounds for fundamental solutions to non-local diffusion
  equations with divergence free drift.
\newblock {\em J. Funct. Anal.}, 264(10):2245--2268, 2013.

\bibitem[Mos64]{Mos64}
J.~Moser.
\newblock A {H}arnack inequality for parabolic differential equations.
\newblock {\em Comm. Pure Appl. Math.}, 17:101--134, 1964.

\bibitem[RC25]{Rui25}
J.~Ruiz-Cases.
\newblock Fractional fast diffusion with initial data a {R}adon measure.
\newblock {\em arXiv:2503.14296}, 2025.

\bibitem[WZ23]{WeZa23}
F.~Weber and R.~Zacher.
\newblock Li-{Y}au inequalities for general non-local diffusion equations via
  reduction to the heat kernel.
\newblock {\em Math. Ann.}, 385(1-2):393--419, 2023.

\bibitem[Wei22]{Wei22}
M.~Weidner.
\newblock Energy methods for nonsymmetric nonlocal operators.
\newblock {\em PhD Thesis (Bielefeld University)}, 2022.

\bibitem[Wid44]{Wid44}
D.~V. Widder.
\newblock Positive temperatures on an infinite rod.
\newblock {\em Trans. Amer. Math. Soc.}, 55:85--95, 1944.

\end{thebibliography}
\end{document}